\documentclass[11pt]{article}

\usepackage{enumerate,xspace}
\usepackage{amsmath,xspace,amssymb,mathrsfs}
\usepackage[dvips]{graphics}

\input xy
\xyoption{all}
\xyoption{2cell}
\UseAllTwocells
\CompileMatrices

\title{Restriction categories III: colimits, partial limits, and extensivity}
\author{J.R.B. Cockett \thanks{Partially supported by NSERC, Canada.} 
\\ Department of Computer Science,
University of Calgary,\\ Calgary, Alberta, T2N 1N4, Canada 
\and Stephen Lack \thanks{Supported by the Australian Research Council.} \\ 
School of Mathematics and Statistics,\\ 
University of Sydney, NSW 2006,\\ Australia}

\renewcommand{\phi}{\varphi}

\newcommand{\X}{\ensuremath{\mathbf X}\xspace}
\newcommand{\Y}{\ensuremath{\mathbf Y}\xspace}
\newcommand{\one}{\ensuremath{\mathbf 1}\xspace}

\newcommand{\C}{{\ensuremath{\mathscr C}}\xspace}
\newcommand{\D}{{\ensuremath{\mathscr D}}\xspace}
\newcommand{\E}{{\ensuremath{\mathscr E}}\xspace}
\newcommand{\M}{{\ensuremath{\mathscr M}}\xspace}
\newcommand{\Q}{{\ensuremath{\mathscr Q}}\xspace}
\newcommand{\R}{{\ensuremath{\mathscr R}}\xspace}
\renewcommand{\S}{{\ensuremath{\mathscr S}}\xspace}
\newcommand{\V}{{\ensuremath{\mathscr V}}\xspace}

\newcommand{\two}{\ensuremath{{\hbox{\textrm 2}\kern-.25em
        \hbox{\vrule height1.5ex width 0.4pt depth -.2ex}}\kern.2em}\xspace}

\newcommand{\Total}{{\sf Total}\xspace}
\newcommand{\rst}[1]{\overline{#1}}
\newcommand{\TotX}{{\sf Total(\X)}\xspace}
\newcommand{\FamC}{{\sf Fam(\ensuremath{\mathscr C})}\xspace}
\newcommand{\FamX}{{\sf Fam(\ensuremath{\mathbf X})}\xspace}
\newcommand{\rid}{{\rm RId}\xspace}
\newcommand{\Dpar}{{\ensuremath{{\mathscr D}_{+1}}}\xspace}
\newcommand{\ExX}{{\sf Ex(\X)}\xspace}
\newcommand{\sLat}{{\sf sLat}\xspace}
\newcommand{\sLatFib}{{\sf sLatFib}\xspace}
\newcommand{\dsLat}{{\sf dsLat}\xspace}
\newcommand{\dsLatFib}{{\sf dsLatFib}\xspace}

\newcommand{\Cat}{{\sf Cat}\xspace}
\newcommand{\rCat}{{\sf rCat}\xspace}
\newcommand{\rCatl}{{\sf rCatl}\xspace}
\newcommand{\erCat}{{\sf erCat}\xspace}
\newcommand{\ParC}{{\ensuremath{\Par_\M(\C)}}\xspace}

\newcommand{\Par}{{\sf Par}\xspace}
\newcommand{\Set}{{\sf Set}\xspace}
\newcommand{\Ab}{{\sf Ab}\xspace}
\newcommand{\CRng}{{\sf CRng}\xspace}
\newcommand{\CRngp}{\ensuremath{\CRng_\t}\xspace}

\newcommand{\Dist}{{\sf Dist}\xspace}
\newcommand{\Distcl}{\ensuremath{\Dist_{\textnormal{cl}}}\xspace}
\newcommand{\Extpr}{\ensuremath{{\sf Ext}_{\textnormal{pr}}}\xspace}
\newcommand{\Copy}{{\sf Copy}\xspace}
\newcommand{\Kcr}{\ensuremath{K_{\textnormal{cr}}}\xspace}
\newcommand{\<}{\langle}
\renewcommand{\>}{\rangle}

\newcommand{\op}{\ensuremath{^{\textnormal{op}}}}

\newcommand{\ot}{\otimes}         
\renewcommand{\t}{\times}         

\newtheorem{theorem}{Theorem}[section]    
\newtheorem{corollary}[theorem]{Corollary}   
\newtheorem{lemma}[theorem]{Lemma}   
\newtheorem{preremark}[theorem]{Remark}   
\newtheorem{preexample}[theorem]{Example}   
\newtheorem{prequestion}[theorem]{Question}   
\newtheorem{proposition}[theorem]{Proposition}
\newtheorem{predefinition}[theorem]{Theorem}

\newenvironment{remark}{\begin{preremark}\rm}{\end{preremark}}
\newenvironment{example}{\begin{preexample}\rm}{\end{preexample}}
\newenvironment{question}{\begin{prequestion}\rm}{\end{prequestion}}

\newcommand{\bm}{\begin{bmatrix}}
\newcommand{\fm}{\end{bmatrix}}

\newcommand{\proof}{\noindent{\sc Proof:}\xspace}
\def\endproof{~\hfill$\Box$\vskip 10pt}

\begin{document}

\maketitle
\begin{abstract}
A restriction category is an abstract formulation for a category of 
partial maps, defined in terms of certain specified idempotents called
the restriction idempotents. All categories of partial maps are restriction 
categories; conversely, a restriction category is a category of partial
maps if and only if the restriction idempotents split. Restriction categories 
facilitate reasoning about partial maps as they have a purely algebraic 
formulation.

In this paper we consider colimits and limits in restriction categories. 
As the notion of restriction category is not self-dual, we should not expect
colimits and limits in restriction categories to behave in the same manner. 
The notion of colimit in the restriction context is quite straightforward,
but limits are more delicate. The suitable notion of limit turns out to
be a kind of lax limit, satisfying certain extra properties.

Of particular interest is the behaviour of the coproduct both by itself 
and with respect to partial products.  We explore various conditions under 
which the coproducts are ``extensive'' in the sense that the total category 
(of the related partial map category) becomes an extensive category.  When 
partial limits are present, they become ordinary limits in the total category. 
Thus, when the coproducts are extensive we obtain as the total 
category a lextensive category.  This provides, in particular, a description 
of the extensive completion of a distributive category.
\end{abstract}


\section{Introduction}


In a category \C with a suitable class \M of monomorphisms, one can define
a category \ParC, of partial maps in \C whose domain of definition lies in
\M. The resulting category has further structure which determines, among
other things, the extent of the partiality involved. This is necessary as
an abstract category can arise as a category of partial maps in more than
one way. For example, any
category \X can be regarded as the category of partial maps in \X where 
the class \M consists only of the isomorphisms (the ``total subobjects'');
thus if $\X=\ParC$, we have $\ParC=\Par_{\textsf Iso}(\X)$.

To describe this extra structure a variety of techniques have been employed.
We recall below four possible approaches to capturing this further 
structure,  indicating 
 how the ``trivial case'' where \M is just the isomorphisms can be
identified.

\begin{enumerate}[A.]
\item Given partial maps $f,g:A\to B$, we define $f\le g$ if $g$ is defined
whenever $f$ is, and they then agree. This makes $\ParC$ into a bicategory,
and is the approach taken by Carboni in \cite{carboni}; it is also closely
related to Freyd's notion of allegory \cite{cats-alligators}. The trivial
case is characterized by the fact that the partial order is discrete, in the
sense that $f\le g$ only if $f=g$.
\item If the category \C of total maps has finite products, then this
induces a symmetric monoidal structure on $\ParC$, given on objects by
the product in \C. The trivial case can be characterized by the fact
that this symmetric monoidal structure on the category of partial maps
is in fact cartesian (that is, given by the categorical product). This
approach was taken by Robinson and Rosolini \cite{robinson-rosolini}
and by Curien and Obtulowicz \cite{curien}.
\item If \C has a strict initial object, and the unique map out of the
initial object is in \M, then $\ParC$ has zero maps, 
given by the ``nowhere defined'' partial maps. These were fundamental
in the approach of di~Paola and Heller \cite{dipaola-heller}. The presence
of these zero maps means that only when the category itself is trivial 
can the partiality be trivial.
\item To every partial map $f:A\to B$ we can associate the partial
map $\rst{f}:A\to A$ which is defined whenever $f$ is, in which it
acts as the identity. This operation is taken as fundamental in the
notion of restriction category studied in
the earlier instalments \cite{resti,restii} of this sequence of papers
and again here. The maps of the form $\rst{f}$ are always 
idempotents, and are called restriction idempotents.
This time the trivial case is characterized by the fact that 
the restriction idempotents are just the identity maps. 
\end{enumerate}

The assignment of $\rst{f}:A\to A$ to $f:A\to B$ mentioned above satisfies
four axioms:
\begin{enumerate}[{[\bf R.1]}]
\item $f\rst{f}=f$ for all $f:A\to B$;
\item $\rst{f}\rst{g}=\rst{g}\rst{f}$ for all $f:A\to B$ and $g:A\to C$;
\item $\rst{g\rst{f}}=\rst{g}\rst{f}$ for all $f:A\to B$ and $g:A\to C$;
\item $\rst{g}f=f\rst{gf}$ for all $f:A\to B$ and $g:B\to C$.
\end{enumerate}

A key property 
of restriction categories, not shared by the axiomatics of 
\cite{dipaola-heller, robinson-rosolini,curien,carboni},
is that any full subcategory of a restriction category has an induced
restriction structure; in fact the restriction categories are precisely 
the full subcategories of categories of partial maps. Conversely, a 
restriction category is a category of partial maps if and only if 
the restriction idempotents split. Restriction categories 
facilitate reasoning about partial maps as they have a purely algebraic 
formulation, which does not involve having any structure on the types.

In this paper we consider the structure on a restriction category arising
from limits and colimits on the category of total maps.
As the notion of restriction category is not self-dual, we should not expect
colimits and limits in restriction categories to behave in the same manner. 
The notion of coproduct in the restriction context is quite straightforward:
a restriction category with {\em restriction coproducts} is just a
cocartesian object in the 2-category \rCat of restriction categories. 
This means that the diagonal $\X\to\X\t\X$ and the 
unique map $\X\to\one$ to the terminal restriction category both have left
adjoints in the 2-category \rCat. This is described in more concrete 
terms in Section~\ref{sect:coproducts}; it means that the category \X
has coproducts which satisfy certain conditions involving the restriction
structure.

On the other hand a {\em cartesian} object in \rCat necessarily 
has a trivial restriction structure. The suitable notion of 
a restriction category with restriction products turns out to be a cartesian 
object in a 2-category \rCatl with the same objects and 1-cells as \rCat, but
with a certain type of ``lax natural transformation'' as 2-cells.
This time the underlying category of a restriction category with restriction
products does not in general have products, 
although the category of total maps does so; a concrete description is 
given in Section~\ref{sect:products}. The resulting structure turns out 
to be equivalent to the p-categories of Robinson and Rosolini 
\cite{robinson-rosolini}, and indeed to other different formulations by a 
variety of authors.
A restriction category can also have products which are entirely independent 
of the restriction structure.  The presence of such products does have 
the slightly surprising effect of ensuring that the lattices of restriction 
idempotents have finite joins over which the meets distribute.

More generally, the suitable notion of limit turns out to be a certain
type of lax limit, and we briefly explore these in Section~\ref{sect:limits}. 
Once again, restriction limits in a restriction category become ordinary 
limits in the category of total maps.

Of particular interest is the behaviour of the coproduct both by itself 
and with respect to partial products.  We explore in 
Section~\ref{sect:extensive} various conditions under 
which the coproducts are ``extensive'' in the sense that the total category 
(of the related partial map category) becomes an extensive category.  When 
partial limits are present, they become ordinary limits in the total category. 
Thus, when the coproducts are extensive we obtain as the total 
category a lextensive category.  This provides, in particular, an 
alternative description of the extensive completion of a distributive category
to that given in \cite{extcomp}. This is described in Section~\ref{sect:Dext}.
But what is the importance of being extensive?  Section~\ref{sect:coproducts}
answers this question for partial map categories very concretely: extensivity 
means that there is a ``calculus of matrices.''  This is critical 
to understanding and manipulating the maps in these settings.

\subsubsection*{Notation}
The identity morphism on an object $A$ is denoted by $A$ or
$1_A$. We write $\<f|g\>:A+B\to C$ for the morphism induced by $f:A\to C$
and $g:B\to C$ and $\<f,g\>:A\to B\t C$ for the morphism induced
by $f:A\to B$ and $g:A\to C$. We also write 
$\<f_\lambda\>:\sum_{\lambda\in\Lambda}A_\lambda\to B$ for the morphism
induced by a $\Lambda$-indexed family of morphisms $f_\lambda:A_\lambda\to B$.
Our notation for coproduct injections
is more flexible: sometimes we write $i$ and $j$ for the two injections
of a binary coproduct, and sometimes we use $i$ with a suitable subscript.
We write $\tau:A+B\to B+A$ for the canonical isomorphism.
The projections of a product are usually denoted by $\pi$ with
a suitable subscript.


\section{Coproducts and matrices}
\label{sect:coproducts}


It is well known that in the category of sets and binary relations
the disjoint union (of a finite family of sets) serves both as coproduct
and product, so that there is a ``calculus of matrices'': see
\cite{matrices}, for example.  In this section we consider the extent
to which this can be adapted to deal not with relations but with partial
functions. We then consider when such a calculus is available in an 
abstract category of partial maps, or restriction category.

Given finite families 
$(A_\lambda)_{\lambda\in \Lambda}$ and $(B_\kappa)_{\kappa\in K}$ of sets,
and a partial function $f:\sum_\lambda A_\lambda\to \sum_\kappa B_\kappa$, we 
may define a partial function $f_{\lambda\kappa}:A_\lambda\to B_\kappa$ for 
each $\lambda\in \Lambda$ and $\kappa\in K$,
by declaring $f_{\lambda\kappa}(x)$ to be defined if and only if $f(x)$ is
defined and lies in $B_\kappa$, in which case $f_{\lambda\kappa}(x)=f(x)$.
Conversely, a matrix 
$(f_{\lambda\kappa})_{\lambda\in \Lambda,\kappa\in K}$, with 
$f_{\lambda\kappa}$ a partial function from $A_\lambda$ to $B_\kappa$ for each 
$\lambda$ and $\kappa$, determines a relation $f$ from 
$\sum_\lambda A_\lambda$ to $\sum_\kappa B_\kappa$, where if 
$x\in A_\lambda$ and $y\in B_\kappa$ we have $f(x)=y$ if and only if
$f_{\lambda\kappa}(x)=y$. The relation $f$ is in fact a partial function
precisely when, for each $\lambda\in\Lambda$, if $f_{\lambda\kappa}(x)$
and $f_{\lambda\kappa'}(x)$ are both defined then $\kappa=\kappa'$: in
other words, if for each $x$ and $\lambda$, there is at most one $\kappa$
for which $f_{\lambda\kappa}(x)$ is defined.

Not only can we represent partial functions by matrices, we can represent
composition of partial functions by matrix multiplication, in the following
sense. If $f:\sum_\lambda A_\lambda\to\sum_\kappa B_\kappa$ and 
$g:\sum_\kappa B_\kappa\to \sum_k C_k$
are partial functions with matrices 
$(f_{\lambda\kappa})_{\lambda\in \Lambda, \kappa\in K}$ and 
$(g_{\kappa\mu})_{\kappa\in K,\mu\in M}$, then the matrix of $gf$ is 
$(\vee_\kappa g_{\kappa\mu}f_{\lambda\kappa})_{\lambda\in \Lambda,\mu\in M}$, 
where $\vee_\kappa g_{\kappa\mu}f_{\lambda\kappa}$ is the partial function 
$h:A_\lambda\to C_k$ with $h(x)=g_{\kappa\mu}f_{\lambda\kappa}(x)$ if the 
right hand side is defined for some (necessarily unique) $\kappa$, and 
undefined otherwise.

If $f$ is defined by $t:E\to\sum_\kappa B_\kappa$ with domain 
$m:E\to\sum_\lambda A_\lambda$, then $f_{\lambda\kappa}$ can be computed as a 
pullback, as in 
$$
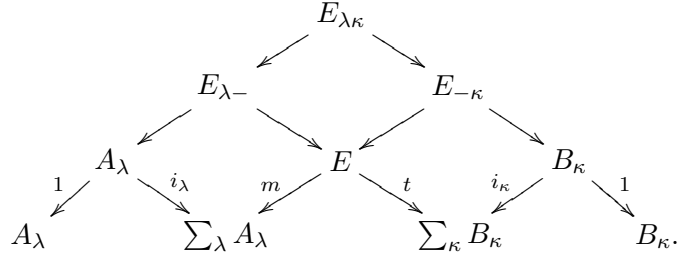
\xymatrix @R1pc @C1pc { 
&&& E_{\lambda\kappa} \ar[dl] \ar[dr] \\
&& E_{\lambda-} \ar[dr] \ar[dl] && E_{-\kappa} \ar[dl] \ar[dr] \\
& A_\lambda \ar[dl]_1 \ar[dr]^{i_\lambda} && E \ar[dl]_m \ar[dr]^t && 
B_\kappa \ar[dl]_{i_\kappa} \ar[dr]^1 \\
A_\lambda && {}\sum_\lambda A_\lambda && {}\sum_\kappa B_\kappa && B_\kappa. 
}$$
In effect we are composing $f$ with the injection 
$i_\lambda:A_\lambda\to\sum_\lambda A_\lambda$, seen as a
total partial map, and the partial map 
$i^*_\kappa:\sum_\kappa B_\kappa\to B_\kappa$ which is defined
as the identity on $B_\kappa$ and is undefined elsewhere. More 
abstractly, $i^*_\kappa$ is the (unique) map satisfying 
$i^*_\kappa i_\kappa=1$ and $i_\kappa i^*_\kappa=\rst{i^*_\kappa}$.
(We shall say that $i^*_\kappa$ is the restriction retraction of 
$i_\kappa$.)

We can recover $f$ from the $f_{\lambda\kappa}$ as the composite
$$\xymatrix{
{\sum_\lambda A_\lambda \ar[r]^{\sum_\lambda h_\lambda} } &
{\sum_{\lambda\kappa} A_\lambda\ar[r]^{\sum_{\lambda\kappa}f_{\lambda\kappa}}}&
{\sum_{\lambda\kappa} B_\kappa \ar[r]^{\sum_\kappa \nabla} } &
{\sum_\kappa B_\kappa} }$$
where $h_\lambda:A_\lambda\to\sum_\kappa A_\lambda$ is defined by
$h_\lambda(x)=(x,\kappa)$ if $f_{\lambda\kappa}(x)$ is defined for
some (necessarily unique) $\kappa$, and undefined otherwise. Once
again, there is also a more abstract characterization of $h_\lambda$:
it is the unique map satisfying $h'_\lambda h_\lambda=\rst{h_\lambda}$
and $h_\lambda h'_\lambda=\rst{h'_\lambda}$, where 
$h'_\lambda:\sum_\kappa A_\lambda\to A_\lambda$ is 
$\<\rst{f_{\lambda\kappa}}\>_{\kappa\in K}$. (We shall say that 
$h_\lambda$ is the restriction inverse of $h'_\lambda$.)

What structure does a restriction category \X need in order to 
support such a calculus of matrices? Obviously \X must have
finite coproducts, and the coproduct injections must have 
restriction retractions. Also, given a morphism 
$f:A\to\sum_\kappa B_\kappa$, the map 
$\<\rst{i^*_\kappa f}\>_{\kappa\in K}:\sum_\kappa A\to A$ must have
a restriction inverse. This sets up a pair of functions between
\begin{itemize}
\item the set of morphisms from $\sum_\lambda A_\lambda$ to 
$\sum_\kappa B_\kappa$, and
\item the set of matrices 
$(f_{\lambda\kappa}:A_\lambda\to B_\kappa)_{\lambda\in\Lambda,\kappa\in K}$
with the property that for each $\lambda$, the map
$\<i^*_\kappa \rst{fi_\lambda}\>_{\kappa\in K}:\sum_\kappa A_\lambda\to A_\lambda$
has a restriction inverse.
\end{itemize}
Finally we need these functions to be mutually inverse and to respect
composition. We shall investigate when this occurs in the remainder of
Section~\ref{sect:coproducts}.


\subsection{Restriction coproducts}


In the previous section we saw that for a restriction category \X with
coproducts to admit a calculus of matrices, it is necessary
that the coproduct injections be restriction monics, and so in particular
be total. In this section we examine the situation in which the coproduct
injections are total. 

\begin{lemma}
Let \X be a restriction category with coproducts, and suppose that
the injections of every binary coproduct $A+B$ are total. Then
\begin{enumerate}[(i)]
\item the unique arrow $z_A:0\to A$ is total for every object $A$;
\item the codiagonal $\nabla:A+A\to A$ is total for every object $A$;
\item $\rst{f+g}=\rst{f}+\rst{g}$ for all arrows $f$ and $g$.
\end{enumerate}
\end{lemma}

\proof
To prove ($iii$), let $f:A\to A'$ and $g:B\to B'$ and write 
$i:A\to A+B$, $j:B\to A+B$, $i':A'+B'$, and $j':B'\to A'+B'$ for
the injections. Then 
$\rst{(f+g)}i=i\rst{(f+g)i}=i\rst{i'f}=i\rst{f}$ since $i'$ is
total, and similarly $\rst{(f+g)}j=j\rst{g}$; thus $\rst{f+g}=\rst{f}+\rst{g}$.
The proof of ($ii$) is similar, while ($i$) follows immediately from the
fact that $z_A$ is an injection of the coproduct $A+0$.
\endproof

We say that such a restriction category has {\em restriction coproducts}.
A more abstract point of view is that such a restriction category 
is just a {\em cocartesian object} in the 2-category \rCat. Recall
that an ordinary category \C has binary coproducts if and only if 
the diagonal functor $\C\to\C\t\C$ has a left adjoint. More generally
an object \C of a 2-category with finite products is said to be 
cocartesian if the diagonal $\C\to\C\t\C$ has a left adjoint in the
2-category. Thus a cocartesian restriction category is a restriction
category \X for which the diagonal
restriction functor $\X\to\X\t\X$ and the unique restriction functor
$\X\to\one$ to the terminal restriction category both have left adjoints
in the 2-category \rCat. Any 2-functor takes adjunctions to adjunctions,
and a finite-product-preserving 2-functor takes cocartesian objects
to cocartesian objects. For instance, there is a 2-functor 
$\Total:\rCat\to\Cat$
which sends a restriction category \X to its category of total maps,
and clearly \Total preserves finite products. Then again, there is
a 2-functor $K_r:\rCat\to\rCat$ which sends a restriction category \X
to the restriction category $K_r(\X)$ obtained by splitting the 
restriction idempotents of \X.

\begin{proposition}\label{prop:cocart} 
If \X is a restriction category with restriction coproducts then
\TotX and $\Total(K_r(\X))$ have coproducts. If $F:\X\to\Y$ is
a coproduct-preserving restriction functor between restriction
categories with restriction coproducts, then $\Total(F):\X\to\Y$
and $\Total(K_r(F)):\Total(K_r(\X))\to\Total(K_r(\Y))$ preserve coproducts.
\end{proposition}

\proof
The 2-functors \Total and $\Total(K_r)$ send cocartesian objects
to cocartesian objects, and so send restriction categories with
restriction coproducts to categories with coproducts.

Similarly, if $F$ preserves coproducts then it commutes with 
the left adjoints $\one\to\X$  and $\X\t\X\to\X$, and so 
$\Total(F)$ commutes with the induced left adjoints $1\to\TotX$ and
$\TotX\t\TotX\to\TotX$; that is, $\Total(F)$ preserves coproducts.
The case of $\Total(K_r(F))$ is entirely analogous.
\endproof

This proposition has a converse when the restriction category is
classified. Recall \cite{resti} that an arrow $r:A\to B$ in a restriction
category is said to be a {\em restriction retraction} if there is an 
arrow $i:B\to A$ with $ri=1$ and $\rst{r}=ir$; such an $i$ is unique.
Recall further \cite{restii} that a restriction category \X is
{\em classified} if the inclusion $\TotX\to\X$ has a right adjoint $R$,
and for each object $A$ the counit $\epsilon_A:RA\to A$ is a restriction
retraction. The promised converse is now:

\begin{proposition}
If \X is a classified restriction category and \TotX has coproducts, 
then \X has restriction coproducts. An arbitrary functor $F:\X\to\C$
preserves coproducts if and only if its restriction $\TotX\to\C$ 
to the total maps preserves coproducts. In particular, for a restriction
category \Y with restriction coproducts, a restriction
functor $F:\X\to\Y$ preserves coproducts if and only if 
$\Total(F):\TotX\to\Total(\Y)$ does so. 
\end{proposition}

\proof
Since \X is classified, the inclusion $\TotX\to\X$ is a left adjoint, and 
so preserves all existing colimits. Since it is also bijective on objects, 
\X has coproducts if \TotX does so; and the injections are
clearly total. 

Since the inclusion $I:\TotX\to\X$ is bijective on objects, a functor
$G:\TotX\to\C$ preserves coproducts if and only if $GI$ does so. Since
$FI$ is just the composite of $\Total(F)$ and the inclusion $\Total(\Y)\to\Y$,
it follows that $F$ preserves coproducts.
\endproof

We also have:

\begin{proposition}
If \X is a restriction category with coproducts and a zero object,
then \X has restriction coproducts.
\end{proposition}

\proof
If \X has a zero object then the injection $i:A\to A+B$
has a retraction $\<1|0\>:A+B\to A$, and so is monic; but
monomorphisms are always total.
\endproof

\begin{example}
If \D is a distributive category, then the endofunctor $+1$ of \D
has a well-known monad structure, and the Kleisli category \Dpar 
of this monad has a restriction structure described in Example~7
of Section~2.1.3 of \cite{resti}. Since \D has coproducts
and the left adjoint $I:\D\to\Dpar$ is bijective on objects, \Dpar
has coproducts; the injections are in the image of $I$ and so
total. Thus \Dpar has restriction coproducts.
\end{example}

It is well known (see \cite{CLW,cockett} for example) that the free 
completion under (finite) coproducts of a category \C can be formed as the 
category \FamC of finite families of objects of \C. Explicitly, an 
object of \FamC is a finite family $(A_\lambda)_{\lambda\in\Lambda}$ of 
objects of \C, and a morphism from $(A_\lambda)_{\lambda\in\Lambda}$ to 
$(B_\kappa)_{\kappa\in K}$ consists
of a function $\phi:\Lambda\to K$ and a family 
$(f_\lambda:A_\lambda\to B_{\phi\lambda})_{\lambda\in I}$
of morphisms in \C. The universal property of \FamC is
expressed in terms of the fully faithful functor $J:\C\to\FamC$ sending
an object of \C to the corresponding singleton family. 

The observation we wish to make here is:

\begin{remark} If \X is a restriction
category then \FamX has a canonical restriction structure, with
$\rst{(\phi,f)} =\left(1_\Lambda,(\rst{f_\lambda})_{\lambda\in\Lambda}\right)$.
Then $J:\X\to\FamX$ is clearly a restriction functor. Furthermore, \X has 
restriction coproducts if and only if $J:\X\to\FamX$ has a left adjoint in 
\rCat. A purely formal consequence is that \FamX is the free restriction 
category with restriction coproducts on \X; we leave the precise formulation 
of the universal property to the reader. Another straightforward 
observation is that the restriction category \FamX is classified
whenever \X is so.  
\end{remark}



\subsection{Restriction zero objects}


To begin with, we allow \X to be an arbitrary restriction category.
Given arrows $f:A\to B$ and $g:B\to A$ in \X,
recall \cite{resti} that $g$ is {\em restriction inverse} to $f$ (and $f$ 
to $g$) if $gf=\rst{f}$ and $fg=\rst{g}$. A restriction inverse is unique if
it exists. In the special case where $f$ is total, we have 
$gf=\rst{f}=1$; then $f$ is said to be a {\em restriction monic} and
$g$ its {\em restriction retraction}, and we often write $f^*$ for $f$.

We say that a zero object $0$ in a restriction category is a 
{\em restriction zero} if for every object $A$ the zero map
$0_{AA}:A\to A$ is a restriction idempotent; that is, $\rst{0_{AA}}=0_{AA}$.

\begin{lemma}\label{lemma:zero}
For a restriction category \X, the following are equivalent:
\begin{enumerate}[{\em ($i$)}]
\item \X has a restriction zero;
\item \X has an initial object $0$ and a terminal object $1$, and each
 $z_A:0\to A$ is a restriction monic;
\item \X has a terminal object $1$ and each $t_A:A\to 1$ is a 
restriction retraction.
\end{enumerate}
\end{lemma}

\proof
$(i)\Rightarrow(ii)$.
If $0$ is a restriction zero then it is both initial and terminal,
and for any object $A$ there is a unique $z_A:0\to A$ and a unique
$t_A:A\to 0$. Clearly $t_Az_A=1$, since $0$ is initial, while
$z_At_A=0_{AA}=\rst{0_{AA}}=\rst{z_At_A}$.

$(ii)\Rightarrow(iii)$.
Let $z^*_1:1\to 0$ be the restriction retraction of $z_1:0\to 1$,
and $t_A:A\to 1$ be the unique map. Then $t_Az_Az^*_1=1$, since $1$ is 
terminal; and we must show that $z_Az^*_1t_A$ is a restriction idempotent.
Now $t_A=z_1z^*_A$, since $1$ is terminal, and so
$z_Az^*_1t_A=z_Az^*_1z_1z^*_A=z_Az^*_A$, which is indeed a restriction 
idempotent.

$(iii)\Rightarrow(i)$.
For each object $A$, choose $s_A:1\to A$ satisfying ($t_As_A=1$ and) 
$\rst{s_At_A}=s_At_A$.
Then $s_0:1\to 0$ is inverse to $z_1:0\to 1$, and so
$0$ is a zero object. Finally $\rst{0_{AA}}=\rst{s_At_A}=s_At_A=0_{AA}$.
\endproof

We now suppose once again that \X has coproducts.

\begin{proposition}
Let \X be a restriction category with coproducts, in which the
coproduct injections are restriction monics. Then the initial object 0
is a restriction zero if and only if the maps $i^*:A+B\to A$ are
natural in $B$; they are always natural in $A$.
\end{proposition}

\proof The $i^*$ can be seen as $1_A+z^*_B:A+B\to A+0$, which are
clearly natural in $A$, and will be natural in $B$ if and only if
the $z^*_B:B\to 0$ are so. But this will be the case if and only if
$0$ is not just initial but also terminal, and now the result follows
by Lemma~\ref{lemma:zero}.
\endproof

We now observe that in order to have a calculus of matrices, the 
category \X must have a restriction zero object. We have 
already seen that the coproduct injections must be restriction
monics, and so in particular that $z_A:0\to A$ must be one.
To deal with empty coproducts,
every map $f:A\to 0$ should be representable as an ``empty matrix'',
which clearly means that there can be at most one such map. Thus in 
this case $0$ is not just an initial object but a zero object (that is,
an initial and a terminal object.) By Lemma~\ref{lemma:zero} it follows
that the initial object is a restriction zero.

\begin{example}
If \D is a distributive category, then the initial object of \D
is a restriction zero in \Dpar. To see this, observe that the left 
adjoint $I:\D\to\Dpar$ 
preserves colimits, so $0$ is initial in \Dpar. For every object 
$A$, there is a unique arrow $A\to 0+1=1$ in \D, and so $0$ is 
also terminal in \Dpar. The zero map $0_{AA}:A\to A$ in \Dpar
is 
$$\xymatrix{A \ar[r]^{!} & 1 \ar[r]^-{i_2} & A+1}$$
and its restriction is
$$\xymatrix{A \ar[r]^-{\<1,i_2!\>} & A\t(A+1) \ar[r]^-{\delta^{-1}} &
A\t A+A \ar[r]^-{\pi_1+!} & A+1.}$$
The fact that these two maps agree is an easy exercise in
distributive categories.
\end{example}

\begin{lemma}
If \X is a restriction category with restriction coproducts
and a restriction zero, then:
\begin{enumerate}[{\em ($i$)}]
\item  each coproduct injection $i:A\to A+B$ is a restriction monic,
       with restriction retraction $i^*:A+B\to A$ equal to $\<1|0\>:A+B\to A$,
       so that the restriction idempotent $ii^*$ is $1+0:A+B\to A+B$;
\item  if $f:C\to A+B$ is total, and the restriction idempotent   
       $\rst{i^*f}$ splits, then the section $k:C_A\to C$
       of the splitting is the pullback in \TotX of the injection 
       $i:A\to A+B$ along $f$;
\item  the natural transformations in \TotX whose components are
       the coproduct injections are cartesian.
\end{enumerate}
\end{lemma}

\proof
($i$) We can regard $i$ as $1_A+z_B$. Then $(1_A+z^*_B)(1_A+z_B)=1$,
while $(1_A+z_B)(1_A+z^*_B)=1_A+0_{BB}=\rst{1_A}+\rst{0_{BB}}=
\rst{1_A+0_{BB}}$.

($ii$) Suppose that $k:C_A\to C$ and $k^*:C\to C_A$ provide the 
splitting, so that $k^*k=1$ and $kk^*=\rst{ii^*f}$.
Then $ii^*fkk^*=(1+0)f\rst{(1+0)f}=(1+0)f=\rst{1+0}f=
f\rst{(1+0)f}=fkk^*$, and so $ii^*fk=fk$. We claim that the commutative
square
$$\xymatrix{C_A \ar[r]^{k} \ar[d]_{i^*fk} & C \ar[d]^f \\ A \ar[r]_-i & A+B}$$
is in fact a pullback in \TotX. Since $i$ and $k$ are monic, it will
suffice to show that a total map $u:D\to C$ factorizes through $k$ 
if $fu$ factorizes through $i$. But if $fu$ factorizes through $i$
then $ii^*fu=fu$, and now
$kk^*u=\rst{ii^*f}u=u\rst{ii^*fu}=u\rst{fu}=u$.

($iii$) We are to show that the square
$$\xymatrix{ A \ar[r]^-{i} \ar[d]_{f} & A + B \ar[d]^{f + g} \\
             A' \ar[r]_-{i'}  & A' + B' }$$
is a pullback in \TotX. Since 
$\rst{i'^*(f+g)}=\rst{fi^*}=\rst{i^*}=ii^*$, the result follows by 
part ($ii$).
\endproof

As we saw above, in the category of sets and relations a coproduct
$A+B$ is also a product, but in the case of sets and partial functions
this is no longer the case. We now describe the trace which remains
of this product structure. In a restriction
category \X with restriction coproducts and a restriction zero, we have 
a functor $+:\X\t\X\to\X$, and natural transformations $i^*:A+B\to A$
and $j^*:A+B\to B$. If there were a natural diagonal $\Delta:A\to A+A$
satisfying the triangle equations, this would exhibit $A+B$ as the 
product of $A$ and $B$. Although there is not such a $\Delta$, we shall
see that there are various maps which ``try''
to be the diagonal; we shall call them decisions.


\subsection{The calculus of matrices}


In this section we consider a restriction category \X with restriction 
coproducts and a restriction zero. The main aim of this section is to 
establish, under further conditions still to be determined, a bijection 
between arrows $f:\sum_\lambda A_\lambda\to\sum_\kappa B_\kappa$ and matrices
$(f_{\lambda\kappa})$ with the property that for each $\lambda$
the map $(\rst{f_{\lambda\kappa}})_\kappa:\sum_\kappa A_\lambda\to A_\lambda$
has a restriction inverse $h_\lambda$. This bijection should send 
$f$ to $(i^*_\kappa fi_\lambda)_{\lambda,\kappa}$ and 
$(f_{\lambda\kappa})_{\lambda,\kappa}$ to the composite
$$\xymatrix @C3pc {
{ \sum_\lambda A_\lambda \ar[r]^-{\sum_\lambda h_\lambda} } &
{ \sum_{\lambda\kappa} A_\lambda \ar[r]^-{\sum_{\lambda\kappa} f_{\lambda\kappa}
               } } &
{ \sum_{\lambda\kappa} B_\kappa \ar[r]^-{\nabla} } & {\sum_\kappa B_\kappa}. }$$

The universal property of the coproduct $\sum_\lambda A_\lambda$ 
reduces this to the case where $\Lambda$ is a singleton.
Thus we are to establish a bijection between the set of morphisms 
$f:A\to\sum_\kappa B_\kappa$ and the set of those 
$K$-tuples $(f_\kappa:A\to B_\kappa)$ for which
$(\rst{f_\kappa})_\kappa:\sum_\kappa A\to A$ has a restriction inverse $h$.
For any $f:A\to\sum_\kappa B_\kappa$ the
induced map $\<\rst{i^*_\kappa f}\>_\kappa:\sum_\kappa A\to A$ will clearly
need to have a restriction inverse $h$. Moreover, $\rst{h}$ will have 
to be $\rst{f}$. For if $h$ is restriction inverse to 
$\<\rst{i^*_\kappa}\>_\kappa$ then
$$\rst{h}=\<\rst{i^*_\kappa f}\>_\kappa h=
\nabla\sum_\kappa\rst{i^*_\kappa f}h=\nabla\rst{\sum_\kappa i^*_\kappa f}h=
\nabla h\rst{(\sum_\kappa i^*_\kappa f)h}$$
but for our bijection we need $(\sum_\kappa i^*_\kappa f)h=f$, so that
$\rst{h}=\nabla h\rst{f}$. But then 
$\rst{h}=\rst{\nabla h\rst{f}}=\rst{h\rst{f}}=\rst{h}\,\rst{f}=
\rst{h}\,\rst{(\sum_\kappa i^*_\kappa f)h}=\rst{(\sum_\kappa i^*_\kappa f)h}=
\rst{f}$, as claimed. 

If $\<\rst{i^*_\kappa f}\>_\kappa$ does have a restriction inverse $h$
and $\rst{h}$ is $\rst{f}$, then we write $[f]$ for $h$, and call it a
{\em decision} for $f$ or {\em $f$-decision}, for reasons which will 
become clearer below. 

\begin{proposition}\label{prop:decision}
An arrow $h:A\to\sum_\kappa A$ is the decision of $f:A\to\sum_\kappa B_\kappa$
if and only if $\nabla h=\rst{f}$ and the square
$$\xymatrix @C3pc {
A \ar[r]^{h} \ar[d]_{f} & {\sum_\kappa A \ar[d]^{\sum_\kappa f}} \\
{\sum_\kappa B_\kappa \ar[r]_{\sum_\kappa i_\kappa} } & 
{\sum_{\kappa,\kappa'\in K} B_\kappa} }$$
commutes.
\end{proposition}

We defer to the next section the proof of the proposition. Observe,
however, that it helps to explain the name ``decision''.
Since $\rst{[f]}=\rst{\nabla[f]}=\rst{f}$, the decision $[f]$ is
defined whenever $h$ is, and the effect
of $[f]$ is ``to send an element $a\in A$ to the element
in the component of $\sum_\kappa A$ corresponding to the component
of $f(a)\in\sum_\kappa B_\kappa$''.

\begin{theorem}\label{thm:matrices}
Let \X be a restriction category with restriction coproducts and 
a restriction zero, in which every map $f:A\to\sum_\kappa B_\kappa$
has a decision. Then there is a bijection between the set of all maps
$f:\sum_\lambda A_\lambda\to\sum_\kappa B_\kappa$ and the set of those matrices
$(f_{\lambda\kappa}:A_\lambda\to B_\kappa)_{\lambda,\kappa}$ for
which $\<\rst{f_{\lambda\kappa}}\>_{\kappa}:\sum_\kappa A_\lambda\to A_\lambda$
has a restriction inverse for every $\lambda$. The bijection sends
$f$ to the matrix $(i^*_\kappa fi_\lambda)_{\lambda,\kappa}$.
\end{theorem}

\proof
Write $\Phi$ for the function computing the matrix of a map
$f:\sum_\lambda A_\lambda\to\sum_\kappa B_\kappa$, and $\Psi$ 
for the purported inverse, which sends $(f_{\lambda\kappa})_{\lambda,\kappa}$ 
to the composite
$$\xymatrix{
{ \sum_\lambda A_\lambda \ar[r]^{\sum_\lambda h_\lambda} } &
{ \sum_{\lambda\kappa} A_\lambda \ar[rr]^{\sum_{\lambda\kappa}
   f_{\lambda\kappa}} } &&
{ \sum_{\lambda\kappa}B_\kappa } \ar[r]^{\nabla} &
{ \sum_\kappa B_\kappa } }$$
where $h_\lambda$ is restriction inverse to 
$\<\rst{f_{\lambda\kappa}}\>_{\kappa}:\sum_\kappa A_\lambda\to A_\lambda$.

Starting with $f:\sum_\lambda A_\lambda\to \sum_\kappa B_\kappa$
we get the matrix $\Phi(f)=(i^*_\kappa fi_\lambda:A_\lambda\to B_\kappa)$; and
then $\Psi(\Phi(f))$ is the composite
$$\xymatrix @C3pc {
{ \sum_\lambda A_\lambda \ar[r]^{\sum_\lambda [fi_\lambda]} } &
{ \sum_{\lambda\kappa} A_\lambda 
  \ar[r]^{\sum_{\lambda\kappa}i^*_\kappa fi_\lambda} } &
{ \sum_{\lambda\kappa} B_\kappa \ar[r]^\nabla } & 
{ \sum_\kappa B_\kappa. } }$$
To see that this is just $f$, observe that in the diagram
$$\xymatrix{
& {\sum_{\lambda\kappa}A_\lambda \ar[r]^-{\sum_{\lambda\kappa}i_\lambda } } &
{ \sum_{\lambda\kappa\lambda'}A_{\lambda'} \ar[r]^{\sum_{\lambda\kappa}f} } &
{ \sum_{\lambda\kappa\kappa'}B_{\kappa'} 
  \ar[dr]^{\sum_{\lambda\kappa}i^*_\kappa} } \\
{ \sum_\lambda A_\lambda \ar[ur]^{\sum_\lambda[fi_\lambda]} \ar[dr]_1
  \ar[r]^{\sum_\lambda i_\lambda} } &
{ \sum_{\lambda\lambda'}A_\lambda \ar[r]^{\sum_\lambda f} \ar[d]^\nabla } &
{ \sum_{\lambda\kappa} B_\kappa \ar[ur]^{\sum_{\lambda\kappa}i_\kappa}
  \ar[rr]_1 } &&
{ \sum_{\lambda\kappa}B_\kappa \ar[d]^\nabla } \\
& {\sum_\lambda A_\lambda \ar[rrr]_f } &&& {\sum_\kappa B_\kappa } }$$
the large upper parallelogram commutes by Proposition~\ref{prop:decision},
the upper triangle commutes since $i^*_\kappa i_\kappa=1$, the lower
triangle by one of the triangle equations, and the large lower rectangle
by naturality of $\nabla$. Thus the entire diagram commutes and
$\Psi(\Phi(f))=f$.

Suppose on the other hand that we are given 
$f_{\lambda\kappa}:A_\lambda\to B_\kappa$ for each 
$\lambda\in\Lambda$ and $\kappa\in K$, and that 
$\<\rst{f_{\lambda\kappa}}\>_\kappa:\sum_\kappa A_\lambda\to A_\lambda$
has a restriction inverse $h_\lambda$ for each $\lambda\in\Lambda$.
Then $\Phi\Psi$ sends the matrix $\<f_{\lambda\kappa}\>_{\lambda,\kappa}$
to the composite 
$$\xymatrix{
{ A_\lambda \ar[r]^-{i_\lambda} } &
{ \sum_\lambda A_\lambda \ar[r]^{\sum_\lambda h_\lambda} } &
{\sum_{\lambda\kappa}A_\lambda \ar[rr]^{\sum_{\lambda\kappa}f_{\lambda\kappa}}}
 && {\sum_{\lambda\kappa} B_\kappa \ar[r]^\nabla } & 
{ \sum_\kappa B_\kappa \ar[r]^{i^*_\kappa} } & B_\kappa }$$
which, by the naturality of $i^*_\lambda$ and the definition of 
$\nabla:\sum_\lambda\kappa B_\kappa\to B_\kappa$, is just
$$\xymatrix{
{ A_\lambda \ar[r]^-{h_\lambda} } &
{ \sum_\kappa A_\lambda \ar[r]^{\sum_\kappa f_{\lambda\kappa}} } &
{ \sum_\kappa B_\kappa \ar[r]^-{i^*_\kappa} } & B_\kappa. }$$
Naturality of $i^*_\kappa$ gives 
$i^*_\kappa(\sum_\kappa f_{\lambda\kappa})=f_{\lambda\kappa}i^*_\kappa$,
thus we must show that 
$f_{\lambda\kappa}i^*_\kappa h_\lambda=f_{\lambda\kappa}$.

Now $i^*_\kappa$ is restriction inverse to $i_\kappa$, and
$h_\lambda$ is restriction inverse to $\<\rst{f_{\lambda\kappa}}\>_\kappa$,
so $i^*_\kappa h_\lambda$ is restriction inverse to 
$\<\rst{f_{\lambda\kappa}}\>_\kappa i_\kappa$, which is just
$\rst{f_{\lambda\kappa}}$. But restriction idempotents are their
own restriction inverses, so $i^*_\kappa h_\lambda=\rst{f_{\lambda\kappa}}$.
Thus $f_{\lambda\kappa}i^*_\kappa h_\lambda=
f_{\lambda\kappa}\rst{f_{\lambda\kappa}}=f_{\lambda\kappa}$, and so
$\Phi\Psi$ is indeed the identity, and the bijection is established.
\endproof

We end this section by showing how to ``multiply'' matrices:

\begin{proposition}
Under the hypotheses of Theorem~\ref{thm:matrices}, if 
$f:\sum_\lambda A_\lambda\to\sum_\kappa B_\kappa$ has matrix
$(f_{\lambda\kappa})_{\lambda,\kappa}$, and 
$g:\sum_\kappa B_\kappa\to\sum_\mu C_\mu$ has matrix
$(g_{\kappa\mu})_{\kappa,\mu}$, then the composite $gf$ has
matrix
$(\vee_\kappa g_{\kappa\mu}f_{\lambda\kappa})_{\lambda,\mu}$,
where $\vee_\kappa g_{\kappa\mu}f_{\lambda\kappa}:A_\lambda\to C_\mu$
is given by 
$$\xymatrix{
A_\lambda \ar[r]^-{h_\lambda} & 
{\sum_\kappa A_\lambda \ar[rr]^{\<g_{\kappa\mu}f_{\lambda\kappa}\>_\kappa} }&&
C_\mu }$$
and $h_\lambda$ is the restriction inverse of 
$\<\rst{f_{\lambda\kappa}}\>_\kappa:\sum_\kappa A_\lambda\to A_\lambda$.
\end{proposition}

\proof
We must show that
$$i^*_\mu gfi_\lambda=\<g_{\kappa\mu}f_{\lambda\kappa}\>_\kappa h_\lambda.$$ 
By the theorem
$fi_\lambda=\sum_\kappa(i^*_\kappa fi_\lambda)h_\lambda$, so
$i^*_\mu gfi_\lambda = i^*_\mu g\sum_\kappa(i^*_\kappa fi_\lambda)h_\lambda
= \<i^*_\mu gi_\kappa i^*_\kappa fi_\lambda\>_\kappa h_\lambda
= \<g_{\kappa\mu}f_{\lambda\kappa}\>_\kappa h_\lambda$
as required.
\endproof


\subsection{Decisions}


In this section we further explore decisions in a restriction
category \X with restriction coproducts and restriction zero; the 
main goal is to prove Proposition~\ref{prop:decision}. Recall that
$h:A\to\sum_\kappa A$ is the decision of $f:A\to\sum_\kappa B_\kappa$
if it is restriction inverse to 
$\<\rst{i^*_\kappa f}\>_\kappa:\sum_\kappa A\to A$ and $\rst{h}=\rst{f}$.
We say that $h:A\to\sum_\kappa A$ is a decision if it is the 
decision of some map $f:A\to\sum_\kappa B_\kappa$.

\begin{example}\label{ex:decisions}~
\begin{enumerate}[($i$)]
\item If $K$ is a singleton, so that we have a single map $f:A\to B$,
a decision for $f$ is a map $h:A\to A$ which is restriction inverse 
to $\rst{f}$: this is just $\rst{f}$ itself.
\item If $K$ is empty, so that $f$ is the unique map $A\to 0$, a decision for
$f$ is a map $h:A\to 0$ which is restriction inverse to the unique
map $z_A:0\to A$: then $f=h=z^*_A$.
\item Let $f$ be a coproduct injection 
$i_\lambda:A_\lambda\to\sum_\lambda A_\lambda$. Then 
$\<\rst{i^*_\kappa i_\lambda}\>_\kappa:\sum_\kappa A_\lambda\to A_\lambda$
is $i^*_\lambda$, which has restriction inverse $i_\lambda$. Thus
$i_\lambda$ is its own decision.
\end{enumerate}
\end{example}

\begin{proposition}
For a map $h:A\to\sum_{\kappa\in K} A$ the following are equivalent:
\begin{enumerate}[{\em ($i$)}]
\item $h$ is its own decision;
\item $h$ is a decision;
\item $h$ has a restriction inverse $g:\sum_\kappa A\to A$ and
$gi_\kappa:A\to A$ is a restriction idempotent for each $\kappa$.
\end{enumerate}
\end{proposition}

\proof
The downward implications are trivial; we must show that given
restriction idempotents $e_\kappa:A\to A$ for each $\kappa\in K$,
if $\<e_\kappa\>_\kappa:\sum_\kappa A\to A$ has a restriction
inverse $h$ then $h$ is its own decision. 

Since $h$ is restriction inverse to $\<e_\kappa\>_\kappa$ and
$i^*_\kappa$ is restriction inverse to $i_\kappa$,
we see that $i^*_\kappa h$ is restriction inverse to 
$\<e_\kappa\>_\kappa i_\kappa$; but the latter is just $e_\kappa$
which is its own restriction inverse. Thus $i^*_\kappa h=e_\kappa$.
and so $\rst{i^*_\kappa h}=e_\kappa$. But then $h$ is restriction
inverse to $\<\rst{i^*_\kappa h}\>_\kappa$, which is just to say
that $h$ is its own decision.
\endproof

The next result says that we can ``conjugate'' decisions by
restriction inverses:

\begin{corollary}\label{cor:conjugate}
If $h:A\to\sum_\kappa A$ is a decision, and $f:A\to B$ a map
with restriction inverse $g:B\to A$, then 
$$\xymatrix{
B \ar[r]^-g & A \ar[r]^-h & {\sum_\kappa A \ar[r]^{\sum_\kappa f} } & 
{\sum_\kappa B} }$$
is a decision and $\rst{(\sum_\kappa g)hf}=\rst{hf}$.
\end{corollary}

\proof
Since $h$ is a decision it is restriction inverse to 
$\<\rst{i^*_\kappa h}\>_\kappa:\sum_\kappa A\to A$. Since
$g$ is restriction inverse to $f$, and $\sum_\kappa f$ is 
restriction inverse to $\sum_\kappa g$, also 
$(\sum_\kappa f)hg$ is restriction inverse to 
$f\<\rst{i^*_\kappa h}\>_\kappa(\sum_\kappa g)$.
Now
$f\<\rst{i^*_\kappa h}\>_\kappa(\sum_\kappa g)i_\kappa=
f\<\rst{i^*_\kappa h}\>_\kappa i_\kappa g=
f\rst{i^*_\kappa h}g=fg\rst{i^*_\kappa hg}=\rst{g}\rst{i^*_\kappa hg}$
which is a restriction idempotent, thus
$(\sum_\kappa f)hg$ is a decision by the Proposition.

Finally, 
$\rst{(\sum_\kappa g)hf}=\rst{\nabla(\sum_\kappa g)hf}=
\rst{g\nabla hf}=\rst{g\rst{h}f}=\rst{gf\rst{hf}}=\rst{\rst{f}\,\rst{hf}}=
\rst{f}\,\rst{hf}=\rst{hf}$.
\endproof

\begin{corollary}\label{cor:decision-on-sum}
If $h:\sum_\lambda A_\lambda\to\sum_\kappa\sum_\lambda A_\lambda$ is
a decision then so is $k_\lambda=(\sum_\kappa i^*_\lambda)hi_\lambda$
for each $\lambda$, and $h$ is the composite
$$\xymatrix{
{\sum_\lambda A_\lambda \ar[rr]^-{\sum_\lambda k_\lambda} } &&
{\sum_\lambda\sum_\kappa A_\lambda \ar[r]^{\sigma} } &
{\sum_\kappa\sum_\lambda A_\lambda} }$$
where $\sigma$ is the canonical isomorphism.
\end{corollary}

\proof
The fact that $k_\lambda$ is a decision is immediate from the previous 
corollary. On the other hand 
$\sigma(\sum_\lambda k_\lambda)i_\lambda=
\sigma i_\lambda k_\lambda=(\sum_\kappa i_\lambda)k_\lambda=
(\sum_\kappa i_\lambda)(\sum_\kappa i^*_\lambda)hi_\lambda=
(\sum_\kappa \rst{i^*_\kappa})hi_\lambda=
hi_\lambda\rst{(\sum_\kappa i^*_\kappa)hi_\lambda}=
hi_\lambda\rst{hi_\lambda}=
hi_\lambda$ 
for each $\lambda$, where the penultimate step uses the previous
corollary. Thus $\sigma(\sum_\lambda k_\lambda)=h$ as claimed.
\endproof

We are now ready to prove Proposition~\ref{prop:decision}.
We shall make frequent use of the naturality of $i^*_\kappa$:

\vskip\baselineskip
\noindent{\sc Proof of Proposition~\ref{prop:decision}:}
First we simplify the condition for $h:A\to\sum_\kappa A$ to be the decision of
$f:A\to\sum_\kappa B_\kappa$. This will be the case if 
$h\<\rst{i^*_\kappa f}\>_\kappa=\rst{\<\rst{i^*_\kappa f}\>_\kappa}$
and $\<\rst{i^*_\kappa f}\>_\kappa h=\rst{h}$. Now
$\<\rst{i^*_\kappa f}\>_\kappa=\nabla(\sum_\kappa\rst{i^*_\kappa f})=
\nabla\rst{\sum_\kappa(i^*_\kappa f)}$ and so the first condition 
becomes 
$$h\rst{i^*_\kappa f}=i_\kappa\rst{i^*_\kappa f}.$$ 

Suppose that $(\sum_\kappa f)h=(\sum_\kappa i_\kappa)f$ and
$\nabla h=\rst{f}$. Then $\rst{h}=\rst{\nabla h}=\rst{f}$.
Now $(\sum_\kappa\rst{f})h=\rst{\sum_\kappa f}h=h\rst{(\sum_\kappa f)h}=
h\rst{(\sum_\kappa i_\kappa)h}=h\rst{h}=h$, and so
\begin{equation*}
\rst{i^*_\kappa h} = \rst{i^*_\kappa(\sum_\kappa\rst{f})h}
                   = \rst{\rst{f}i^*_\kappa h} 
                   = \rst{fi^*_\kappa h} 
                   = \rst{i^*_\kappa(\sum_\kappa f)h} 
                   = \rst{i^*_\kappa(\sum_\kappa i_\kappa)f}
                   = \rst{i_\kappa i^*_\kappa f} 
                   = \rst{i^*_\kappa f}
\end{equation*}
but now
\begin{equation*}
h\rst{i^*_\kappa f} = h\rst{i^*_\kappa h} 
                    = \rst{i^*_\kappa}h 
                    = i_\kappa i^*_\kappa h 
                    = i_\kappa\nabla i_\kappa i^*_\kappa h 
                    = i_\kappa\nabla\rst{i^*_\kappa} h 
                    = i_\kappa\nabla h\rst{i^*_\kappa h} 
                    = i_\kappa\rst{h}\,\rst{i^*_\kappa h} 
                    = i_\kappa\rst{i^*_\kappa h}
\end{equation*}
giving the first condition. As for the second
\begin{multline*}
\<\rst{i^*_\kappa f}\>_\kappa h 
                    = \nabla(\sum_\kappa\rst{i^*_\kappa f})h 
                    = \nabla\rst{\sum_\kappa(i^*_\kappa f)} h 
                    = \nabla h\rst{(\sum_\kappa i^*_\kappa f)h} \\
                    = \rst{h} \rst{(\sum_\kappa i^*_\kappa)(\sum_\kappa f)h}
            = \rst{h} \rst{(\sum_\kappa i^*_\kappa)(\sum_\kappa i_\kappa)f} 
            = \rst{h}\, \rst{f} 
            = \rst{h}
\end{multline*}
and so $h$ is the decision of $f$.

Suppose conversely that $h$ is the decision of $f$. Then
\begin{equation*}
\rst{\sum_\kappa i^*_\kappa f}h 
                = h\rst{(\sum_\kappa i^*_\kappa f)h} 
                = h\rst{\rst{\sum_\kappa i^*_\kappa f}h}
                = h\rst{\<\rst{i^*_\kappa f}\>_\kappa h} 
                = h\rst{h} = h
\end{equation*}
and so $\nabla h=\nabla\rst{\sum_\kappa i^*_\kappa f}h=
\nabla(\sum_\kappa\rst{i^*_\kappa f})h=\<\rst{i^*_\kappa f}\>_\kappa h=\rst{h}=
\rst{f}$. On the other hand
\begin{equation*}
(\sum_\kappa f)h  = (\sum_\kappa f)\rst{\sum_\kappa i^*_\kappa f}h
          = \sum_\kappa( f\rst{i^*_\kappa f})h 
          = \sum_\kappa(\rst{i^*_\kappa}f)h
          = (\sum_\kappa \rst{i^*_\kappa})(\sum_\kappa f)h
\end{equation*}
and so 
\begin{multline*}
(\sum_\kappa i_\kappa)f 
          = (\sum_\kappa i_\kappa)f\rst{h}
          = (\sum_\kappa i_\kappa)f\nabla\rst{\sum_\kappa i^*_\kappa f}h
          = (\sum_\kappa i_\kappa)\nabla(\sum_\kappa f)
              \rst{(\sum_\kappa i^*_\kappa)(\sum_\kappa f)}h \\
   = (\sum_\kappa i_\kappa)\nabla\rst{\sum_\kappa i^*_\kappa}(\sum_\kappa f)h
   = (\sum_\kappa i_\kappa)\nabla(\sum_\kappa i_\kappa i^*_\kappa f)h
   = (\sum_\kappa i_\kappa)(\sum_\kappa i^*_\kappa f)h
   = (\sum_\kappa \rst{i^*_\kappa})(\sum_\kappa f)h 
   = (\sum_\kappa f)h.
\end{multline*}
\endproof

We saw in Example~\ref{ex:decisions} that a decision for
$f:A\to\sum_{\kappa\in K}B_\kappa$ always exists if $K$ is 
empty or a singleton. We end this section by proving that
all decisions exist provided that binary ones do.

\begin{proposition}
A restriction category \X with restriction coproducts and
a restriction zero has all decisions provided that it
has a decision for each $f:A\to B+C$.
\end{proposition}

\proof
Let $f:A\to\sum_{\kappa\in K}B_\kappa$ be given, where $K$
is a finite set of cardinality greater than 2. Choose $\lambda\in K$,
and regard $\sum_{\kappa\in K}B_\kappa$ as the coproduct of 
$B_\lambda$ and $\sum_{\kappa\neq\lambda}B_\kappa$ with injections
$i$ and $j$. By assumption, 
$f:A\to B_\lambda+(\sum_{\kappa\neq\lambda}B_\kappa)$ has 
a decision $h_\lambda:A\to A+A$. Suppose by way of inductive hypothesis that 
$j^*f:A\to\sum_{\kappa\neq\lambda}B_\kappa$ has a decision
$h':A\to\sum_{\kappa\neq\lambda}A$.  We shall show that
$$\xymatrix{
A \ar[r]^-{h_\lambda} & A+A \ar[r]^-{1+h'} & A+\sum_{\kappa\neq\lambda}A =
\sum_{\kappa\in K}A }$$
is a decision for $f:A\to\sum_{\kappa\in K}B_\kappa$.

Commutativity of
$$\xymatrix{
A \ar[r]^{h_\lambda} \ar[d]_f & A+A \ar[rr]^{1+h'} \ar[d]_{f+f} &&
A+\sum_{\kappa\neq\lambda}A \ar[d]_{f+\sum_{\kappa\neq\lambda}f} \ar@{=}[dr] \\
B_\lambda+\sum_{\kappa\neq\lambda}B_\kappa \ar[r]_{i+j} \ar@{=}[dr] &
{ \sum_{\kappa'}B_{\kappa'}+\sum_\kappa B_\kappa 
    \ar[rr]_{1+\sum_{\kappa\neq\lambda}i_\kappa} } &&
{\sum_{\kappa'}B_{\kappa'}+\sum_{\kappa\neq\lambda,\kappa'}B_{\kappa'} 
    \ar@{=}[dr] }
& {\sum_\kappa A \ar[d]^{\sum_\kappa f} } \\
& { \sum_\kappa B_\kappa \ar[rrr]_{\sum_\kappa i_\kappa} } &&& 
{ \sum_{\kappa,\kappa'}B_{\kappa'} } }$$
gives one of the conditions in Proposition~\ref{prop:decision}; it
remains to show that $\nabla(1+h')h_\lambda=\rst{f}$. Commutativity
of
$$\xymatrix{
A \ar[r]^-{h_\lambda} \ar[dr]_{\rst{(1+j^*f)h_\lambda}} &
A+A \ar[r]^-{1+h'} \ar[dr]_{1+\rst{j^*f}} &
A+\sum_{\kappa\neq\lambda}A \ar[d]^{1+\nabla} \\
& A \ar[r]_{h_\lambda} \ar[dr]_f & A+A \ar[d]^\nabla \\
&& A}$$
reduces this to proving that $f\rst{(1+j^*f)h_\lambda}=\rst{f}$.

To do so, first observe that
$\rst{(i^*f+1)(1+j^*f)h_\lambda}=\rst{(i^*f+j^*f)h_\lambda}=
\rst{\<\rst{i^*f}|\rst{j^*f}\>h_\lambda}=\rst{h_\lambda}=\rst{f}$
so that 
$f\rst{(1+j^*f)h_\lambda}=f\rst{f}\,\rst{(1+j^*f)h_\lambda}=
f\rst{(i^*f+1)(1+j^*f)h_\lambda}\,\rst{(1+j^*f)h_\lambda}=
f\rst{(i^*f+1)(1+j^*f)h_\lambda}=f\rst{(i^*f+j^*f)h_\lambda}=
\rst{f}$
as required.
\endproof

Finally, we record the following result which will be needed below:

\begin{proposition}\label{prop:adding-decisions}
If $f:A\to B+C$ and $f':A'\to B'+C'$ have decisions $h$ and $h'$ then 
$(1+\tau+1)(f+f'):A+A'\to(B+B')+(C+C')$ has decision
$(1+\tau+1)(h+h')$.
\end{proposition}

\proof
Let $i:B\to B+C$, $j:C\to B+C$, $i':B'\to B'+C'$, and $j':C'\to B'+C'$
be the various injections. Then the injection $k:B+B'\to B+B'+C+C'$ is
given by $(1+\tau+1)(i+i')$. Similarly, write $l$ for the injection
$(1+\tau+1)(j+j'):C+C'\to B+B'+C+C'$. 

First observe that 
$\rst{(1+\tau+1)(h+h')}=\rst{h+h'}=\rst{h}+\rst{h'}=\rst{f}+\rst{f'}=
\rst{f+f'}=\rst{(1+\tau+1)(f+f')}$. 
Now $h$ is restriction inverse to $\<\rst{i^*f}|\rst{j^*f}\>$ and
$h'$ is restriction inverse to $\<\rst{i'^*f'}|\rst{j'^*f'}\>$, thus
$h+h'$ is restriction inverse to 
$\<\rst{i^*f}|\rst{j^*f}\>+\<\rst{i'^*f'}|\rst{j'^*f'}\>$, and
$(1+\tau+1)(h+h')$ is restriction inverse to 
$\left(\<\rst{i^*f}|\rst{j^*f}\>+\<\rst{i'^*f'}|\rst{j'^*f'}\>\right)
(1+\tau+1)$. But
\begin{align*}
\left(\<\rst{i^*f}|\rst{j^*f}\>+\<\rst{i'^*f'}|\rst{j'^*f'}\>\right)(1+\tau+1)
&= (\nabla+\nabla)(\rst{i^*f}+\rst{j^*f}+\rst{i'^*f'}+\rst{j'^*f'})(1+\tau+1)\\
&= (\nabla+\nabla)(1+\tau+1)(\rst{i^*f}+\rst{i'^*f'}+\rst{j^*f}+\rst{j'^*f'})\\
&= \nabla(\rst{i^*f+i'^*f'}+\rst{j^*f+j'^*f'}) \\
&= \nabla(\rst{(i+i')^*(f+f')}+\rst{(j+j')^*(f+f')}) \\
&= \<\rst{(i+i')^*(f+f')}|\rst{(j+j')^*(f+f')}\>
\end{align*}
so that $(1+\tau+1)(h+h')$ is the decision of $(1+\tau+1)(f+f')$ 
as claimed.
\endproof


\section{Extensive restriction categories}
\label{sect:extensive}



\subsection{Extensivity}


In the previous section we saw that a restriction category \X admits
a calculus of matrices if it has restriction coproducts, a restriction
zero, and decisions. In this section we relate this structure to the 
question of when \TotX and $\Total(K_r(\X))$ are extensive.

\begin{proposition}\label{prop:total-extensive}
If \X is a restriction category with restriction coproducts and a 
restriction zero, then \TotX is extensive if and only 
if, for every total arrow $f:C\to A+B$, the restriction idempotent
$\rst{(1+0)f}$ splits and an $f$-decision exists.
If \X has an object $1$ which is terminal in \TotX, then it 
suffices to consider the case $A=B=1$.
\end{proposition}

\proof
We know that \TotX has coproducts since \X has restriction coproducts,
and we know that the coproduct injections in \TotX are cartesian, since
\X has a restriction zero. Thus \TotX will be extensive if and only if
it has pullbacks along coproduct injections, and coproducts are stable.

Suppose that $\rst{(1+0)f}$ splits for every $f:C\to A+B$, and that an
$f$-decision exists. Let $k:C_A\to C$ and $k^*C\to C_A$ provide the
splitting for $\rst{(1+0)f}$. Let $l:C_B\to C$ and $l^*:C\to C_B$ 
provide the splitting for $\rst{(0+1)f}$, which exists since 
$\rst{(0+1)f}=\rst{\tau(0+1)f}=\rst{(1+0)\tau f}$. We are to show that
$\<k|l\>:C_A+C_B$ is invertible.

The $f$-decision $h:C\to C+C$ is restriction inverse to $\<kk^*|ll^*\>$,
so that $h\<kk^*|ll^*\>=kk^*+ll^*$ and $\<kk^*|ll^*\>h=\rst{h}=\rst{f}=1$.
Thus $h\<k|l\>=k+l$ and so $(k^*+l^*)h\<k|l\>=(k^*+l^*)(k+l)=1$, while
$\<k|l\>(k^*+l^*)h=\<kk^*|ll^*\>h=1$, as required.

Suppose conversely that \TotX is extensive. Then any map $C\to A+B$
has the form $f+g:A'+B'\to A+B$, and now $if=(f+g)i'$ so that $f$ is
total, and similarly $g$ is total. Also 
$\rst{(1+0)(f+g)}=\rst{i^*(f+g)}=\rst{fi'^*}=\rst{i'^*}$
so that $i'$ and $i'^*$ provide a splitting for $\rst{(1+0)(f+g)}$.
Finally the identity $A'+B'\to A'+B'$ is easily seen to be a decision
for $f+g$.

If $1$ is terminal in $\Total(K(\X))$ then it suffices to show
stability of the coproduct $1+1$; see \cite{CLW} or \cite{cockett}.
\endproof

\begin{corollary}\label{cor:total-Kr-extensive}
If \X is a restriction category with restriction coproducts and a 
restriction zero, then $\Total(K_r(\X))$ is extensive if and only 
every arrow $f:C\to A+B$ in \X has a decision map.
If \X has an object $1$ which is terminal in $\Total(K_r(\X))$, then it 
suffices to consider the case $A=B=1$.
\end{corollary}

\proof
Since \X has restriction coproducts, so does $K_r(\X)$, and since
\X has a restriction zero, so does $K_r(\X)$. All restriction idempotents
split in $K_r(\X)$, so by Proposition~\ref{prop:total-extensive},
$\Total(K_r(\X))$ will be extensive if and only if every total
arrow $f:(C,e)\to(A+B,e_1+e_2)$ has a decision. To say that $f$
is total is to say that $\rst{f}=e$. 

A decision $h$ for $f:(C,e)\to(A+B,e_1+e_2)$ is an arrow $h:C\to C+C$ 
in \X satisfying $(\rst{f}+\rst{f})h=h=h\rst{f}$, $\nabla h=\rst{f}$,
and $(f+f)h=(i+j)f$; that is, a decision map for $f$ in \X.
\endproof

In light of the proposition, we say that a restriction category \X
is {\em extensive} if it has restriction coproducts and a restriction
zero, and every map $f:C\to A+B$ has a decision. By the uniqueness
of decisions and the characterization of Proposition~\ref{prop:decision},
the existence of these decisions can be viewed as a
combinator assigning to each $f:C\to A+B$ a map $\<f\>:C\to C+C$
satisfying the decision axioms:
\begin{enumerate}[{[\bf D.1]}]
\item $\nabla\<f\>=\rst{f}$;
\item $(f+f)\<f\>=(i+j)f$.
\end{enumerate}
Thus decision structure is equational, and so can be added freely.
It would be interesting to have a description of the free extensive
restriction category on a restriction category, or the free such on
a mere category.

Of course to say that \X is extensive as a restriction category is
quite different to saying that is extensive as a mere category. In
fact as an extensive restriction category has a zero object it cannot
be an extensive category unless it is the trivial category with a single
object and a single arrow. The connection between extensive restriction
categories and extensive categories is rather 
(see Corollary~\ref{cor:total-Kr-extensive}) that if \X is an extensive
restriction category then $\Total(K_r(\X))$ is an extensive category.

\begin{example}\label{ex:Dex}
If \D is a distributive category, then \Dpar is an extensive
restriction category, and so $\Total(K_r(\Dpar))$ is an
extensive category. We have already seen that \Dpar has
restriction coproducts and a restriction zero, thus we may
apply Corollary~\ref{cor:total-Kr-extensive}. If $f:C\to A+B+1$
is an arrow in \Dpar from $C$ to $A+B$, let $h$ be the composite
$$\xymatrix{
C \ar[r]^-{\<C,f\>} & C\t(A+B+1) \ar[r]^{\delta^{-1}} &
C\t A+C\t B+C \ar[rr]^-{\pi_1+\pi_1+!} && C+C+1.}$$
Verification of the commutativity of the diagrams
$$\xymatrix{
C \ar[r]^-{\<C,f\>} \ar@/_3pc/[drrr]_{\rst{f}} & 
C\t(A+B+1) \ar[dr]_{\delta^{-1}} \ar[r]^{\delta^{-1}} &
C\t A+C\t B+C \ar[r]^-{\pi_1+\pi_1+!} & C+C+1 \ar[d]^{\nabla+1} \\
&& C\t(A+B)+C \ar[u]_{\delta^{-1}+C} \ar[r]^-{\pi_1+!} & C+1}$$
$$\xymatrix{
C \ar[r]^-{\<C,f\>} \ar[d]_{f} & 
C\t(A+B+1) \ar[r]^-{\delta^{-1}} \ar[d]^{f\t(A+B+1)} &
C\t A+C\t B+C \ar[d]^{\pi_1+\pi_1+!} \\
A+B+1 \ar[r]^-{\Delta} \ar@/_3pc/[rrdd]_{i+j+1} & 
(A+B+1)\t(A+B+1) \ar[d]^-{\delta^{-1}} & C+C+1 \ar[d]^{\nabla+1} \\
& *\txt{$(A+B+1)\t A+(A+B+1)\t B$\\ $+(A+B+1)$} \ar[r]_-{\pi_1+\pi_1+!} &
A+B+1+A+B+1+1 \ar[d]^{t} \\
&& A+B+A+B+1}$$
is a straightforward exercise in distributive categories; here 
$t:A+B+1+A+B+1+1\to A+B+A+B+1$ is the composite of the twist map
$A+B+1+A+B+1+1\to A+B+A+B+1+1+1$ and $A+B+A+B+!$. Thus $h$ is
the required decision for $f$.
\end{example}


\subsection{Extensive maps}
\label{sect:extensive-maps}


As well as considering when \TotX or $\Total(K_r(\X))$ is extensive,
we can look at subcategories which are extensive. To this end, we
say that {\em the  map $f:A\to B$ in \X is extensive} if for any
decision $h:B\to B+B$ there is an $hf$-decision $k:A\to A+A$. 

\begin{lemma}\label{lemma:extensive-map}~
\begin{enumerate}[{\em ($i$)}]
\item Restriction isomorphisms are extensive;
\item Restriction idempotents are extensive;
\item Decision maps are extensive;
\item Coproduct injections and codiagonals are extensive.
\end{enumerate}
\end{lemma}

\proof
First we show that restriction isomorphisms are extensive.
If $f:A\to B$ has restriction inverse $g:B\to A$, and $h:B\to B+B$
is a decision, then $(g+g)hf$ is a decision and
$\nabla(g+g)hf=\rst{(g+g)hf}=\rst{hf}$ by Corollary~\ref{cor:conjugate}.
Thus
\begin{multline*}
(hf+hf)(g+g)hf =(hfg+hfg)hf=(h\rst{g}+h\rst{g})hf=(h+h)\rst{(g+g)}hf \\
=(h+h)hf\rst{(g+g)hf}=(h+h)hf\rst{hf}=(h+h)hf=(i+j)hf
\end{multline*}
and so $(g+g)hf$ is an $hf$-decision, and $f$ is extensive.

Every restriction idempotent is restriction inverse to itself, and
is therefore extensive. Similarly, decisions and coproduct injections are 
restriction isomorphisms and therefore extensive. As for the codiagonal, if
$h:A\to A+A$ is a decision, consider the composite
$$\xymatrix{A+A \ar[r]^-{h+h} & A+A+A+A \ar[r]^{1+\tau+1} & A+A+A+A.}$$
On the one hand we have 
$\nabla(1+\tau+1)(h+h)=(\nabla+\nabla)(h+h)=\rst{h}+\rst{h}=
\rst{h+h}=\rst{\nabla(h+h)}=\rst{h\nabla}$, and on the other,
$(h\nabla+h\nabla)(1+\tau+1)(h+h)=(h+h)\nabla(h+h)=(h+h)h\nabla=
(i+j)h\nabla$; thus $(1+\tau+1)(h+h)$ is
an $h\nabla$-decision.
\endproof

\begin{proposition}
Let \X be a restriction category with restriction coproducts and
a restriction zero. Then the extensive maps in \X form a restriction
subcategory \ExX of \X which is closed under finite coproducts,
contains the decisions; and $\Total(K_r(\ExX))$ is extensive.
Furthermore, \ExX is maximal among restriction subcategories
of \X with these properties.
\end{proposition}

\proof
By Lemma~\ref{lemma:extensive-map}, we know that \ExX contains the identities, 
the restriction idempotents, the coproduct injections and the codiagonals. 
Thus it will be a restriction subcategory provided that it is closed under
composition, and it will be closed under finite coproducts provided that the
extensive maps are so. By Lemma~\ref{lemma:extensive-map} once again,
we know that \ExX contains the decisions; while the fact that 
$\Total(K_r(\ExX))$ is extensive and the maximality of \ExX will follow from 
Corollary~\ref{cor:total-Kr-extensive}. Thus we need only show that the
extensive maps are closed under composition and coproducts.

If $f:A\to B$ and $g:B\to C$ are extensive, and $h:C\to C+C$ is
a decision, let $k:B\to B+B$ be an $hg$-decision, and let
$l:A\to A+A$ be a $kf$-decision. Then
$\nabla l=\rst{kf}=\rst{\nabla kf}=\rst{\rst{hg}f}=\rst{hgf}$ and
\begin{align*}
(hgf+hgf)l &= (hg+hg)(f+f)l \\
           &= (hg+hg)(\rst{hg}+\rst{hg})(f+f)l \\
           &= (hg+hg)(\nabla k+\nabla k))(f+f)l \\
           &= (hg\nabla+hg\nabla)(kf+kf)l \\
           &= (hg\nabla+hg\nabla)(i+j)kf \\
           &= (hg+hg)kf \\
           &= (i+j)hgf
\end{align*}
so that $l$ is also an $hgf$-decision.

Finally, let $f:A\to B$ and $f':A'\to B'$ be extensive, and let
$h:B+B'\to B+B'+B+B'$ be a decision. By Corollary~\ref{cor:decision-on-sum},
$h$ can be written as $(1+\tau+1)(k+k')$, where $k:B\to B+B$ and
$k':B'\to B'+B'$ are decisions. Since $f$ and $f'$ are extensive,
$kf$ and $k'f'$ have decisions, and so, by 
Proposition~\ref{prop:adding-decisions},
$(1+\tau+1)(kf+k'f')$ has a decision; but 
$(1+\tau+1)(kf+k'f')=h(f+f')$, and so we have proved that $f+f'$ is
extensive.
\endproof

Clearly $\Total(K_r(\X))$ is extensive if and only if $\ExX=\X$; that is,
if every map is extensive. Note, however,
that the construction \ExX is not functorial in \X.


\section{Limits in restriction categories}
\label{sect:products}


We saw in Section~\ref{sect:coproducts} that cocartesian objects
in \rCat give a good notion of restriction category with coproducts.
We now turn to products, and the first thing to observe is that
cartesian objects in \rCat are {\em not} a good notion.

If \X is a restriction category, and the unique restriction functor
$!:\X\to\one$ has a right adjoint in \rCat, then \X has a terminal object
$1$, and for each object $A$, the unique map $t_A:A\to 1$ is total. 
But if $f:A\to B$ is any map, then 
$\rst{f}=\rst{\rst{t_B}f}=\rst{t_Bf}=\rst{t_A}=1$,
and so $f$ is total. Thus $!:\X\to\one$ can have a right adjoint in \rCat
only if the restriction structure on \X is trivial. 

The situation for binary products is much the same. Suppose that
$\Delta:\X\to\X\t\X$ has a right adjoint in \rCat. Explicitly,
this means that \X has binary products as a mere category, the
diagonal and projections are total, and $\rst{f\t g}=\rst{f}\t\rst{g}$
for any maps $f$ and $g$. Let $f:A\to B$ be any map, and let 
$p,q:A\t A\to A$ be the projections. Then 
$\rst{1_A\t f}=\rst{\rst{p}(1_A\t f)}=\rst{p(1_A\t f)}=\rst{p}=1$,
and so $1_A\t f:A\t A\to A\t B$ is total; and now
$\rst{f}=\rst{f}q\Delta=q(1_A\t\rst{f})\Delta=q\rst{1_A\t f}\Delta=q\Delta=1$,
so $f$ is total. Thus once again $\Delta:\X\to\X\t\X$ can have a 
right adjoint in \rCat only if the restriction structure on \X is trivial.

We shall now look at other possible notions of products in restriction
categories; and, more generally, limits.


\subsection{Cartesian objects in \rCatl}


One possible approach to the unsatisfactory nature of cartesian
objects in \rCat is to change the 2-category \rCat. In \cite{restii}
we defined a 2-category \rCatl with the same objects and arrows
as \rCat, namely the restriction categories and restriction functors,
but with a larger class of 2-cells. For restriction functors
$F,G:\X\to\Y$, a 2-cell in \rCatl from $F$ to $G$ consists of a
total map $\alpha_X:FX\to GX$ in \Y for each object $X$ of \X,
such that for each $f:X\to Y$ in \X, the diagram
$$\xymatrix{FX \ar[rr]^{Ff} \ar[d]_{F\rst{f}} && FY \ar[d]^{\alpha_Y} \\
FX \ar[r]_{\alpha_X} & GX \ar[r]_{Gf} & GY }$$
commutes. The reason for the name \rCatl is that if one thinks
of a restriction category \X as a 2-category (where there is a
2-cell $f\le g$ if and only if $f=g\rst{f}$) and restriction functors
as 2-functors, then a 2-cell in \rCatl is precisely a lax natural 
transformation from $F$ to $G$ whose components are total.

We now define a {\em restriction terminal object} in a restriction
category \X to be an object $T$ for which the corresponding restriction
functor $\one\to\X$ is right adjoint in \rCatl to the unique functor
$\X\to\one$. In more explicit terms, this amounts to giving, for each
object $A$ of \X, a total map $t_A:A\to T$, such that $t_T=1_T$ and
for each arrow $f:A\to B$, we have $t_B f=t_A \rst{f}$.

\begin{proposition}
A restriction terminal object in \X is terminal in \TotX. Conversely,
if \X is a classified restriction category, then a terminal object
in \TotX is restriction terminal in \X.
\end{proposition}

\proof 
The first statement follows immediately from the fact that  
$\Total:\rCatl\to\Cat$ is a 2-functor, and so preserves adjunctions; 
alternatively, it is equally easy to verify directly.

The second statement is an instance of \cite[Proposition~3.7]{restii}.
\endproof

This proposition means in particular that there is no ambiguity in
saying ``$T$ is a restriction terminal object'', since the total
maps $t_A:A\to T$ are unique. It also shows that restriction terminal
objects are unique up to a unique isomorphism.

Another point of view on restriction terminal objects may be obtained
by consideration of the functor $\rid:\X\op\to\Set$, defined in \cite{resti}.
This sends an object $A$ to the set of all restriction idempotents on
$A$, and a morphism $f:A\to B$ to the function sending a restriction
idempotent $e:B\to B$ to $\rst{ef}:A\to A$.

\begin{proposition}
A restriction terminal object is precisely a representation of
the functor $\rid:\X\op\to\Set$.
\end{proposition}

\proof
If $T$ is a restriction terminal object, then any $f:A\to T$
determines a restriction idempotent $\rst{f}$ on $A$, while 
any restriction idempotent $e:A\to A$ determines a map $t_Ae:A\to T$.
These processes are inverse, since $\rst{t_Ae}=\rst{e}=e$,
and $t_A\rst{f}=t_T f=1_T f=f$.

On the other hand if $T$ is an object equipped with an isomorphism
$\alpha:\X(-,T)\cong\rid$ then for each $A$ there is a unique
$t_A:A\to T$ with $\alpha_A(t_A)=1_A$. For any $f:A\to B$ we
have $\alpha_A(t_Bf)=\rid(f)\alpha_B(t_B)=\rid(f)(1_B)=\rst{f}$
and $\alpha_A(t_A\rst{f})=\rid(\rst{f})\alpha_A(t_A)=\rid(\rst{f})(1_A)
=\rst{f}$; thus $t_Bf=t_A\rst{f}$, since $\alpha_A$ is invertible.
It remains to show that $t_T=1_T$. Let $e$ be the restriction idempotent
$\alpha_T(1_T)$. Then for any $g:A\to T$ we have
$\alpha_A(g)=\alpha_A(1_Tg)=\rid(g)(\alpha_T(1_T))=\rid(g)(e)=\rst{eg}$.
In particular $\alpha_T(e)=\rst{ee}=e=\alpha_T(1_T)$, so that $e=1_T$;
and now $\alpha_T(t_T)=1_T=e=\alpha_T(1_T)$, so that $t_T=1_T$.
\endproof

Next we turn to the case of a restriction category \X for which
the diagonal restriction functor $\Delta:\X\to\X\t\X$ has a right
adjoint in \rCatl. We then say that \X has {\em binary restriction
products}. Explicitly, this means that there is a restriction 
functor $\X\t\X\to\X$ whose value at an object $(A,B)$ we denote $A\t B$ and
whose value at an arrow $(f,g)$ we denote $f\t g$; and total maps
$\Delta:A\to A\t A$, $p:A\t B\to A$, and $q:A\t B\to B$ satisfying
$$\xymatrix{
& A \ar[dl]_1 \ar[dr]^1 \ar[d]_{\Delta} & && A\t B \ar[d]_{\Delta} \ar[dr]^1\\
A & A\t A \ar[l]^{p} \ar[r]_{q} & A && A\t B\t A\t B \ar[r]_-{p\t q} & A\t B \\
A\t B \ar[d]_p & A\t B \ar[l]_{\rst{f}\t\rst{g}} \ar[r]^{\rst{f}\t\rst{g}}
\ar[dd]^{f\t g} & A\t B \ar[d]^q && A \ar[r]^{\rst{f}} \ar[dd]_f &
A \ar[d]^{\Delta} \\
A \ar[d]_f && B \ar[d]^g &&& A\t A \ar[d]^{f\t f} \\
A' & A'\t B' \ar[l]^p \ar[r]_q & B' && A' \ar[r]_{\Delta} & A'\t A'. }$$
Once again, $\Total:\rCatl\to\Cat$ preserves products and adjunctions,
so that \TotX will have binary products whenever $\X\to\X\t\X$ has
a right adjoint in \rCatl.

It turns out that if $\Delta:\X\to\X\t\X$ does have a right adjoint
in \rCatl, then it automatically satisfies certain further conditions,
as the following proposition shows. In particular, the diagonal
maps $\Delta:A\to A\t A$ are not just lax natural, but natural.

\begin{proposition}
If \X is a restriction category, and $\Delta:\X\to\X\t\X$ has a
right adjoint in \rCatl, then: {\em ($i$)} $\rst{(f\t g)\Delta}=\rst{f}\rst{g}$
for all $f$ and $g$ with the same domain, and {\em ($ii$)} the maps
$\Delta:A\to A\t A$ are natural in $A$.
\end{proposition}

\proof
($i$) Since 
$\rst{(f\t g)\Delta}=\rst{\rst{(f\t g)}\Delta}=\rst{(\rst{f}\t\rst{g})\Delta}$,
it will suffice to show that $\rst{(e\t e')\Delta}=ee'$, for all restriction
idempotents $e$ and $e'$.

First observe that 
$\rst{(e\t e')\Delta}=p\Delta\rst{(e\t e')\Delta}=p\rst{e\t e'}\Delta=
p(e\t e')\Delta$, and similarly $\rst{(e\t e')\Delta}=q(e\t e')\Delta$.
Using lax naturality of $p$, we have
$e\rst{(e\t e')\Delta}=ep(e\t e')\Delta=p(e\t e')\Delta=\rst{(e\t e')\Delta}$,
and using lax naturality of $q$, we have
$e'\rst{(e\t e')\Delta}=\rst{(e\t e')\Delta}$. Thus
$\rst{(e\t e')\Delta}=ee'\rst{(e\t e')\Delta}=\rst{(e\t e')\Delta}ee'=
\rst{(e\t e')\Delta ee'}=\rst{(eee'\t e'ee')\Delta ee'}=
\rst{(ee'\t ee')\Delta ee'}=\rst{\Delta ee'}=\rst{ee'}=ee'$.

($ii$) This follows from ($i$) and lax naturality of $\Delta$, since
$(f\t f)\Delta=(f\t f)\Delta\rst{(f\t f)\Delta}=(f\t f)\Delta\rst{f}\,\rst{f}=
(f\t f)\Delta\rst{f}=\Delta f$.
\endproof

We say that \X has restriction products if it is a cartesian object
in \rCatl; that is, if it has binary restriction products and a 
restriction terminal. If \X and \Y are restriction categories with
restriction products, then a restriction functor $F:\X\to\Y$ is
said to preserve restriction products if it commutes with the
right adjoints $\one\to\X$ and $\X\t\X\to\X$ in \rCatl. This definition
can be made more explicit. If $T$ and $S$ denote the restriction
terminal objects of \X and \Y, then there is a unique total map
$\phi:FT\to S$, and $F$ preserves the restriction terminal object
if and only if $\phi:FT\to S$ is invertible. Similarly, for any
objects $X$ and $Y$ of \X there is a unique total map
$\psi_{X,Y}:F(X\t Y)\to FX\t FY$ commuting with the projections,
and $F$ preserves binary restriction products if and only if each
$\psi_{X,Y}$ is invertible. We now have:

\begin{proposition}\label{prop:cart}
If \X is a restriction category with restriction products then
\TotX and $\Total(K_r(\X))$ have products; if \Y is another 
such restriction category and $F:\X\to\Y$ is a restriction functor
which preserves restriction products then $\Total(F):\Total(\X)\to\Total(\Y)$
and $\Total(K_r(F)):\Total(K_r(\X))\to\Total(K_r(\Y))$ preserve products.
\end{proposition}


\subsection{p-Categories}


Here we recall Robinson and Rosolini's notion of 
{\em p-category} \cite{robinson-rosolini}, in order to
compare it to the various structures considered above.

A {\em p-category} is a category \X equipped with a 
functor $\t:\X\t\X\to\X$, a natural family of maps
$\Delta:A\to A\t A$, and families $p_{A,B}:A\t B\to A$ natural
in $A$, and $q_{A,B}:A\t B\to B$ natural in $B$,
required to make commutative the following diagrams:
$$\xymatrix{
& X \ar[dl]_1 \ar[dr]^1 \ar[d]_{\Delta} & && X\t Y \ar[d]_{\Delta} \ar[dr]^1\\
X & X\t X \ar[l]^{p} \ar[r]_{q} & X && X\t Y\t X\t Y \ar[r]_-{p\t q} & X\t Y}$$
$$\xymatrix{
& X\t(Y\t Z) \ar[dl]_{1\t p} \ar[d]_p \ar[dr]^{1\t q} & &&
& (X\t Y)\t Z \ar[dl]_{p\t 1} \ar[d]_q \ar[dr]^{q\t 1} \\
X\t Y \ar[r]_p & X & X\t Z \ar[l]^p && X\t Z \ar[r]_q & Z & Y\t Z \ar[l]^q }$$
$$\xymatrix @C.8pc {
X\t(Y\t Z) \ar[r]^-{\Delta} \ar[d]_{f\t(g\t h)} & 
(X\t(Y\t Z))\t(X\t(Y\t Z)) \ar[rr]^-{(1\t p)\t q} &&
(X\t Y)\t(Y\t Z) \ar[r]^-{1\t q} & (X\t Y)\t Z \ar[d]^{(f\t g)\t h} \\
X'\t(Y'\t Z') \ar[r]^-{\Delta} & (X'\t(Y'\t Z'))\t(X'\t(Y'\t Z')) 
\ar[rr]^-{(1\t p)\t q}&&(X'\t Y')\t(Y'\t Z') \ar[r]^-{1\t q} &(X'\t Y')\t Z'}$$
$$\xymatrix{
X\t Y \ar[d]_{f\t g} \ar[r]^-{\Delta} & (X\t Y)\t(X\t Y) \ar[r]^-{q\t p} &
Y\t X \ar[d]^{f'\t g'} \\
X'\t Y' \ar[r]^-{\Delta} & (X'\t Y')\t(X'\t Y') \ar[r]^-{q'\t p'} & Y'\t X' }$$
for all arrows $f$, $g$, and $h$.
The last two diagrams provide a natural associativity isomorphism
$\alpha_{X,Y,Z}:X\t(Y\t Z)\to(X\t Y)\t Z$ and a natural symmetry 
isomorphism $\tau_{X,Y}:X\t Y\to Y\t X$.

Given a map $f:X\to X'$, Robinson and Rosolini define $\text{dom}f:X\to X$
to be 
$$\xymatrix{X \ar[r]^-{\Delta} & X\t X \ar[r]^{1\t f} & X\t X' \ar[r]^p & X}$$
and their Proposition~1.4 verifies that this makes \X into a restriction
category.
As Robinson and Rosolini observe (in slightly different terminology), 
although a p-category structure on a category may not be unique, a
p-category structure on a restriction category is. Thus it makes 
sense to ask which restriction categories are p-categories.

\begin{proposition}
A restriction category is a p-category if and only if it has
binary restriction products.
\end{proposition}

\proof
First suppose that \X is a p-category.
It is proved in \cite[Proposition~1.4]{robinson-rosolini} that 
$\t:\X\t\X\to\X$ is a restriction functor, and that each instance
of $p$, $q$, and $\Delta$ is total. The diagonal is natural by
assumption, and the ``triangle equations'' linking $\Delta$ with
$p$ and $q$ hold by assumption. Thus it remains only to check that 
the projections $p$ and $q$ are lax natural. In the case of 
$p$, lax naturality amounts to the equation $p(f\t g)=fp(\rst{f}\t\rst{g})$ for
all arrows $f$ and $g$. Consider first the special case where $f$ is
the identity. In the diagram
$$\xymatrix{
X\t Y \ar[r]^-{1\t\Delta} \ar[dr]_{1\t g} & 
X\t(Y\t Y) \ar[r]^{1\t(1\t g)} \ar[dr]_(0.3){1\t g\t g} &
X\t(Y\t Y') \ar[r]^-{1\t p} \ar[d]^{1\t(g\t 1)} & X\t Y \ar[d]^{1\t g} \\
& X\t Y' \ar[r]^-{1\t\Delta} \ar@/_2pc/[rr]_1 & X\t(Y'\t Y') \ar[r]^-{1\t p} & 
X\t Y' \ar[r]^p & X }$$
the left square commutes by naturality of $\Delta$, the triangle
by functoriality of $\t$, the right square by (one-sided) naturality
of $p$, and the curved region by one of the triangle equations. Thus 
the exterior commutes, which is to say that $p(1\t g)=p(1\t\rst{g})$.
As for the general case, in the diagram
$$\xymatrix{
X\t Y \ar@/_2pc/[dd]_{\rst{f}\t\rst{g}} \ar[r]^{f\t g} \ar[d]^{\rst{f}\t1} & 
X'\t Y' \ar[dr]^p \\
X\t Y \ar[r]^{1\t g} \ar[d]^{1\t\rst{g}} & X\t Y \ar[u]^{f\t1} \ar[d]^p & X'\\
X\t Y \ar[r]_p & X \ar[ur]_f }$$
the left and top regions commute by functoriality of $\t$, the bottom
region commutes by the special case just considered, and the right region
by the one-sided naturality of $p$. Commutativity of the exterior is 
the desired equation $p(f\t g)=fp(\rst{f}\t\rst{g})$.

Lax naturality of $q$ states that $q(f\t g)=gq(\rst{f}\t\rst{g})$; we
leave the verification to the reader.

Now suppose conversely that \X has binary restriction products.
We must show that $p:X\t Y\to X$ is natural in $X$, that $q:X\t Y\to Y$ 
is natural in $Y$, and that $\alpha_{X,Y,Z}$ and $\tau_{X,Y}$ are natural
in all variables. The equations involving only instances of
$p$, $q$, and $\Delta$ all hold because the binary restriction
products are actual products in \TotX.

For naturality of $p$, we use lax naturality of $p$ and 
naturality of $\Delta$
to see that $p(f\t 1)=fp\rst{f\t 1}=fp\rst{fp}\rst{q}=fp$; the
case of $q$ is similar.

As for $\tau$, first observe that
$\rst{(g\t f)\tau}=\rst{(g\t f)(q\t p)\Delta}=\rst{(gq\t fp)\Delta}=
\rst{gq}\rst{fp}=\rst{g\t f}$. Now 
$(g\t f)\tau=(g\t f)\tau\rst{(g\t f)\tau}=(g\t f)\tau\rst{g\t f}=\tau(f\t g)$.

The case of $\alpha$ is similar but more complicated. Since
$$\rst{((f\t g)\t h)\alpha}=\rst{((f\t g)\t h)((1\t p)\t qq)\Delta}=
  \rst{((f\t gp)\t hqq)\Delta}=\rst{f\t gp}\,\rst{hqq}=
  \rst{fp}\rst{gpq}\rst{hqq}$$
and
$$\rst{f\t(g\t h)}=\rst{fp}\,\rst{(g\t h)q}=\rst{fp}\,\rst{\rst{g\t h}q}=
  \rst{fp}\,\rst{\rst{gp}\rst{hq}q}=\rst{fp}\rst{q\rst{gpq}\rst{hqq}}=
  \rst{fp}\rst{gpq}\rst{hqq}$$
we have $\rst{((f\t g)\t h)\alpha}=\rst{f\t(g\t h)}$, and now we deduce
$$((f\t g)\t h)\alpha=((f\t g)\t h)\alpha\rst{((f\t g)\t h)\alpha}=
  ((f\t g)\t h)\alpha\rst{f\t(g\t h)}=\alpha(f\t(g\t h)).$$
\endproof

If \X is a p-category, Robinson and Rosolini define a {\em one-element
object} to be an object $T$ with a family $t_X:X\to T$ of maps in \X
for which $p:X\t T\to X$ is invertible, with inverse
$$\xymatrix{X \ar[r]^-{\Delta} & X\t X \ar[r]^-{1\t t_X} & X\t T.}$$

\begin{proposition}
If \X is a p-category, an object $T$ of \X is a one-element object
if and only if it is a restriction terminal object; the map
$t_X:X\to T$ in the definition of one-element object is the unique
total map from $X$ to $T$.
\end{proposition}

\proof
To say that 
$$\xymatrix{X \ar[r]^-{\Delta} & X\t X \ar[r]^-{1\t t_X} & X\t T \ar[r]^-p & 
X}$$
is the identity is precisely to say that $t_X$ is total. The fact that
the $t_X$ are lax natural and $t_T=1_T$ for a one-element object $T$
is part of \cite[Theorem~3.3]{robinson-rosolini}. 

Conversely, if $T$ is a restriction terminal object, then we have a 
family $t_X:X\to T$ of total maps; it remains to show that
$$\xymatrix{X\t T\ar[r]^-p & X\ar[r]^-{\Delta} & X\t X \ar[r]^-{1\t t_X} &
X\t T}$$
is the identity. But this follows from the fact that restriction products
in \X are genuine products in \TotX.
\endproof

The relationship between p-categories with one-element object
and various other structures is analyzed in some detail in
\cite{robinson-rosolini}. Translating this into our nomenclature,
the restriction categories with restriction products are exactly the
{\em partial cartesian categories} in the sense of Curien and Obtulowicz 
\cite{curien}, and or alternatively the {\em pre-dht-symmetric categories} of
Hoehnke \cite{hoehnke}, and they are a special
case of the {\em bicategories of partial maps} of Carboni \cite{carboni}.
For more details on these correspondences, see \cite{robinson-rosolini}.


\subsection{Categories with products and a restriction}


Before leaving our discussion of products in restriction categories,
it is worth discussing a quite different type of product that sometimes
exists. While restriction products are tensor products which
are actual products in the total map category, it is also possible
that a restriction category could have products in the ordinary sense.
A well-known example of this is provided by the category of sets and 
partial maps. The restriction product of $A$ and $B$ is just their
product $A\t B$ as sets; this is not the categorical product in the
partial map category, which is given by $A + A \t B + B$.  
(This is more generally true for the partial map 
category of any lextensive category \cite{CLW,cockett}, where the 
\M-maps are taken to be the coproduct injections.)
   
Recall that the restriction idempotents associated with a particular 
object $A$ in a restriction category \X form a meet semi-lattice with 
$e_1 \wedge e_2 = e_1e_2$, and greatest element $\top = 1_A$. 
These lattices sit over each object to give the {\em restriction fibration} 
\cite[Section~4]{resti}  $\partial: \R(\X) \to \X$ over all the 
maps: the substitutions preserve the meet but not the greatest element. 
When this fibration is restricted to the total maps one obtains a 
meet-semilattice fibration.

If a category has both a restriction structure and a terminal object $1$ 
(and here we emphasize we do not assume any relation between the two 
structures) then we may consider the restriction $\rst{!_A}:A\to A$ of the 
unique map $!_A:A\to 1$. Then for any $f:A\to B$ we have
$\rst{!_A}\,\rst{f}=\rst{!_A\rst{f}}=\rst{!_A}$
and thus $\rst{!_A}$ must be the least restriction idempotent in the above
ordering: in terms of partial maps, this determines the smallest 
possible domain.  Thus the presence of a terminal object forces each object 
$A$ to have a least element in its lattice of restriction idempotents.  
Furthermore, it is clear that the substitution functors of the fibration 
mentioned above preserve these least elements.

When a category has both a restriction structure and finite products 
then $\rst{\<e,e'\>}$ is the {\em join} of $e$ and $e'$ in the lattice 
of restriction idempotents of an object. To see this, first observe that
$$e\rst{\<e,e'\>e}=p\<e,e'\>\rst{\<e,e'\>}=p\<e,e'\>=e$$
so that $e\le\rst{\<e,e'\>}$, and similarly $e'\le\rst{\<e,e'\>}$. Now
if $d$ is a restriction idempotent and $e,e'\le d$, then
$$\rst{\<e,e'\>}d=\rst{\<e,e'\>d}=\rst{\<ed,e'd\>}=\rst{\<e,e'\>}$$
and so $\rst{\<e,e'\>}\le d$.

This proves that the semilattices of restriction idempotents are lattices. 
In fact they are distributive lattices, since
$$e\wedge(e_1\vee e_2)=\rst{\<e_1,e_2\>}e=\rst{\<e_1,e_2\>e}=
\rst{\<e_1e,e_2e\>}=e_1e\vee e_2e=(e\wedge e_1)\vee(e\wedge e_2)$$
and 
$$e\wedge\bot=\rst{!_A}e=\rst{!_Ae}=\rst{!_A}=\bot.$$

\begin{proposition} 
If \X is a category with a restriction structure and (finite) products 
then the fibration of restriction idempotents
$$\partial: \R(\X) \to \X$$
is a fibred join-semilattice and the fibration 
$$\partial_t: \R_t(\X) \to \Total(\X)$$
is a fibred distributive lattice.
\end{proposition}

\proof
It remains only to check that the inverse image functors preserve
the relevant structure. In \cite[Section 4.1]{resti} it was proved
that binary meets are always preserved, while the top element is
preserved by the total maps. Thus it will suffice to show that
for an arbitrary map $f:X\to Y$, the induced functor 
$\rid(f):\rid(Y)\to\rid(X)$ preserves finite joins. For the bottom
element we have $\rid(f)(\rst{!_Y})=\rst{\rst{!_Y}f}=\rst{!_Yf}=\rst{!_X}$; 
for binary joins we have:
$\rid(f)(e\vee e')=\rid(f)\rst{\<e,e'\>}=\rst{\<e,e'\>f}=\rst{\<ef,e'f\>}=
ef\vee e'f=\rid(f)(e)\vee\rid(f)(e')$.
\endproof

In a split restriction category with products this means that the 
$M$-subobjects in the total category must already have finite joins 
which are preserved by pulling back.  Thus {\em products} in the
restriction category lead to {\em colimits} in the lattices of
\M-subobjects.


\subsection{Restriction limits}
\label{sect:limits}


We have already discussed products and coproducts in restriction
categories. Now we turn briefly to more general notions of limit. 
Once again, these will be analyzed in terms of adjunctions in \rCatl.

Let \X be a restriction category and \C a finite category. We shall define
the restriction limit of a functor $S:\C\to\X$ to be a cone
$p_C:L\to SC$ over $S$ with total components, satisfying
the following universal property. If $q_C:M\to SC$ is a lax cone over
$S$ --- that is, $Sc.q_C=q_{D}\rst{Sc.q_C}$ for any $c:C\to D$ --- then
there is a unique arrow $f:M\to L$ satisfying $p_C f=q_C e$, where
$e$ is the composite of the restriction idempotents $\rst{q_C}$.

It follows immediately from the definition that restriction limits
are unique up to unique isomorphism. Equally immediate is the fact
that if $S$ takes its values in \TotX, then a restriction 
limit of \X is a genuine limit in \TotX.

\begin{example}
The restriction limit of the empty diagram is precisely a restriction
terminal object. The restriction limit of a diagram on the discrete
category with two objects is the restriction product of the corresponding
objects.
\end{example}

The following proposition provides a new example:

\begin{proposition}
The restriction limit of an arrow $f$ is precisely a splitting for the
idempotent $\rst{f}$.
\end{proposition}

\proof
A restriction limit of $f$ amounts to a monomorphism $p:P\to X$ for which
$fp$ is total, having the property that for any arrows 
$q:Q\to X$ and $q':Q\to Y$ satisfying $q'=fq\rst{q'}$,
there is a unique $r:Q\to P$ satisfying $pr=qe$ and $fpr=q'e$,
where $e=\rst{q}\rst{q'}$. In fact 
$\rst{q}\rst{q'}=\rst{q}\rst{fq\rst{q'}}=\rst{q}\rst{fq}\rst{q'}=
\rst{fq}\rst{q'}=\rst{q'}$, and $pr=qe$ implies
$fpr=fqe=fq\rst{q'}=q'=q'\rst{q'}=q'e$, and so the only condition on 
$r$ is that $pr=q\rst{q'}$.

Taking $q=1_X$ and $q'=f$, we obtain a unique $s:X\to P$ satisfying
$ps=\rst{f}$. Taking $q=p$ and $q'=fp$, we obtain a unique $t:P\to P$
satisfying $pt=p\rst{fp}=p$. Since $psp=\rst{f}p=p\rst{fp}=p$, we deduce
by uniqueness of $t$ that $sp=1$. Thus $p$ and $s$ provide a splitting 
for $\rst{f}$.

On the other hand, if $p:P\to X$ and $s:X\to P$ split $\rst{f}$, while
$q$ and $q'$ satisfy $q'=fq\rst{q'}$, then 
$psq\rst{q'}=\rst{f}q\rst{q'}=q\rst{fq}\rst{q'}=q\rst{fq\rst{q'}}=q\rst{q'}$,
and $sq\rst{q'}$ is unique with this property, since $p$ is monic.
Thus $p:P\to X$ exhibits $P$ as the restriction limit of $f$.
\endproof

We now, as promised, analyze these restriction limits in terms
of adjunctions in \rCatl. We continue to suppose that \X is a 
restriction category, and now allow \C to be an arbitrary category,
not necessarily finite. We shall define a restriction category
$\X^\C$ and a restriction functor $\Delta:\X\to\X^\C$, and show
that if \C is finite then this $\Delta$ has a right adjoint if
and only if \X has restriction limits of functors with domain \C.

As a category,
$\X^\C$ consists of functors from \C to \X and lax natural transformations
between them. More explicitly, given functors $F,G:\C\to\X$, an arrow 
$\alpha:F\to G$ in $\X^\C$ consists of an arrow $\alpha_A:FA\to GA$ in \X
for each object $A$ of \C, such that 
$\alpha_B.Ff=Gf.\alpha_A\rst{\alpha_B.Ff}$
for every arrow $f:A\to B$ in \C. Composition is defined pointwise:
$(\beta\alpha)_A=\beta_A\alpha_A$. The restriction structure is also defined 
pointwise:
$\rst{\alpha}:F\to F$ has $\rst{\alpha}_A=\rst{\alpha_A}$. The only thing
to check is that $\rst{\alpha}$ is an fact in arrow of the category. To 
do this, note that $\rst{\alpha_B.Ff}=\rst{Gf.\alpha_A.\rst{\alpha_B.Ff}}=
\rst{Gf.\alpha_A}.\rst{\alpha_B.Ff}=
\rst{\alpha_A}.\rst{Gf.\alpha_A}.\rst{\alpha_B.Ff}=
\rst{\alpha_A}.\rst{\alpha_B.Ff}$, and now
$\rst{\alpha_B}.Ff=Ff.\rst{\alpha_B.Ff}=Ff.\rst{\alpha_A}.\rst{\alpha_B.Ff}=
Ff.\rst{\alpha_A}.\rst{\rst{\alpha_B}.Ff}$. The restriction functor 
$\Delta:\X\to\X^\C$ sends an object $X$ of \X to the functor $\C\to\X$
constant at $X$, and sends an arrow $f:X\to Y$ to the family of 
arrows $X\to Y$, each of which is just $f$.

\begin{proposition}
If $\Delta:\X\to\X^\C$ has a right adjoint in \rCatl, 
then \TotX has ordinary \C-limits.
\end{proposition}

\proof
Applying the 2-functor $\Total:\rCatl\to\Cat$ to the adjunction gives an
adjunction
$$\xymatrix{\Total(\X^\C) \ar@/^1pc/[r] \ar@{}[r]|{\top} &
\TotX \ar@/^1pc/[l]^{\Total(\Delta)} }$$
of categories. The functor $\Total(\Delta)$ lands in the
full subcategory $\TotX^\C$ of $\Total(\X^\C)$ consisting
of the functors $\C\to\X$ landing in \TotX. It follows that
$\Delta:\TotX\to\TotX^\C$ has a right adjoint, and so that \TotX
has \C-limits.
\endproof

In order to compare the two approaches to restriction limits,
we assume that the category \C is finite:

\begin{proposition}
If \C is a finite category and \X a restriction category, then
to give a right adjoint in \rCatl to $\Delta:\X\to\X^\C$ is precisely
to give a restriction limit in \X of each functor $S:\C\to\X$.
\end{proposition}

\proof
Let $R:\X^\C\to\X$ be right adjoint in \rCatl to $\Delta$.
If $S:\C\to\X$ is an object of $\X^\C$, let $L=R(S)$, and let
$p:L\to S$ be the component at $S$ of the counit. Then $p$ is a 
lax cone, and its components $p_C:L\to SC$ are total. But to say
that $p$ is a lax cone is to say, for each $f:C\to D$, that
$p_D=Sf.p_C\rst{p_D}$, and since $p_D$ is total, this means that
$p$ is in fact a cone. 

We now show that the cone $p:L\to S$ is a restriction limit cone.
If $q:M\to S$ is a lax cone; that is, an arrow $\Delta M\to S$ in
$\X^\C$, then let $f:M\to R(S)=L$ be given by $R(q):R\Delta M\to R(S)$
composed with the unit $n:M\to R\Delta M$. In the diagram
$$\xymatrix{
\Delta M \ar[r]^{\Delta n} \ar[d]_{\Delta(\rst{Rq.n})} &
\Delta R\Delta M \ar[d]^{\Delta{\rst{Rq}}} \ar[rr]^{\Delta Rq} &&
\Delta RS \ar[d]^{p} \\
\Delta M \ar[r]^{\Delta n} \ar@/_2pc/[rr]_{1}  &
\Delta R\Delta M \ar[r]^{p} & \Delta M \ar[r]_q & S}$$
the left square commutes by a restriction category axiom, 
the right rectangle by lax naturality of $p$, and the curved
region by one of the triangle equations. Commutativity of the
exterior amounts to the equation
$p_C f=q_C. \rst{Rq.n}=q_C\rst{f}$.

We shall now show that $\rst{f}$ is the composite of the $\rst{q_C}$,
which we henceforth denote $e$. For each $C$ we have
$\rst{f}=\rst{p_C f}=\rst{q_C\rst{f}}=\rst{q_C}\rst{f}$, and so
$\rst{f}=e\rst{f}$. On the other hand
$\rst{q_C} e=e$, so $\rst{q}.\Delta e=\Delta e$, and 
$\rst{Rq.R\Delta e}=\rst{R(\rst{q}.\Delta e)}=\rst{R\Delta e}$;
thus 
$e\rst{f}=e\rst{Rq.n}=e\rst{Rq.n.e}=e\rst{Rq.R\Delta e.n.e}=
e.\rst{R\Delta e.n.e}=e.\rst{n.e}=e$. This proves that $\rst{f}$ is
the composite of the $\rst{q_C}$, and so that $f$ provides the
desired factorization.

Finally we must prove that the factorization $f$ is unique. To do this,
we shall show that an arrow $f:M\to RS$ is determined by the $p_Cf$ and by 
$\rst{f}$; then if $p_Cf=p_Cf'$, we also have
$\rst{f}=\rst{p_Cf}=\rst{p_Cf'}=\rst{f'}$. Now in
$$\xymatrix{
M \ar[r]^{n} & R\Delta M \ar[r]^{R\Delta f} & 
R\Delta RS \ar[r]^{Rp} & RS \\
M \ar[u]^{\rst{f}} \ar[rr]_xf&& RS \ar[u]^n \ar[ur]_1 }$$
the rectangle commutes by lax naturality of $n$, and the triangle
by one of the triangle equations; commutativity of the exterior
confirms that $f$ is determined by $p.\Delta f$ and $\rst{f}$, that
is, by the $p_C f$ and by $\rst{f}$. This completes the construction
of restriction \C-limits in \X.

Suppose conversely that \X has restriction \C-limits. For each
$S:\C\to\X$, define $R(S)$ to be the restriction limit of $S$,
and define the component at $S$ of the counit $\Delta R\to 1$ to
be the restriction limit cone $p:\Delta R(S)\to S$. If
$\sigma:S\to T$ is an arrow in $\X^\C$, then for each object $C$
of \C, let $q_C=\sigma_C p_C:R(S)\to TC$. If $f:C\to D$ is an 
arrow of \C, in the following diagram
$$\xymatrix{
R(S) \ar[r]^{p_C} & SC \ar[r]^{\sigma_C} & TC \ar[dr]^{Tf} \\
R(S) \ar[u]^{\rst{\sigma_D.Sf.p_C}} \ar[r]^{p_C} \ar@/_2pc/[rr]_{p_D} &
SC \ar[r]^{Sf} \ar[u]^{\rst{\sigma_D.Sf}} & SD \ar[r]^{\sigma_D} & TD }$$
the left square commutes by one of the restriction category axioms,
the right rectangle by lax naturality of $\sigma$, and the curved
region by naturality of $p$. Finally 
$\rst{\sigma_C.Sf.p_C}=\rst{\sigma_D.p_D}$ by naturality of $p$ once
again, and so $Tf.\sigma_C.p_C.\rst{\sigma_D.p_D}=\sigma_D.p_D$, that
is, $Tf.q_C.\rst{q_D}=q_D$; and so the $q$ form a lax cone. We now
define $R(\sigma):R(S)\to R(T)$ to be the unique arrow for which
$\rst{R(\sigma)}$ is the composite of the restriction idempotents
$\rst{\sigma_C.p_C}$, and the diagram
$$\xymatrix{
R(S) \ar[d]_{\rst{R(\sigma)}} \ar[rr]^{R(\sigma)} && R(T) \ar[d]^{p_C} \\
R(S) \ar[r]_{p_C} & SC \ar[r]_{\sigma_C} & TC }$$
commutes.

The unit $n:M\to R\Delta M$ is defined to be the unique arrow
satisfying $p_Cn=1$, for each leg $p_C:R\Delta M\to M$ of the 
restriction limit cone of $\Delta M$.

We leave to the reader the various straightforward verifications:
that $R$ is a restriction functor, that the unit and counit are
lax natural, and that the triangle equations hold.
\endproof

\begin{proposition}
A restriction category \X has all (finite) restriction limits if and only 
if \X is split as a restriction category and \TotX has finite limits.
\end{proposition}

\proof
We have already seen that \TotX has \C-limits if \X has restriction
\C-limits; and that restriction idempotents split in \X if \X has
restriction \two-limits. Thus it remains to show that if \X is 
a split restriction category and \TotX has finite limits, then
\X has restriction limits.

Let $S:\C\to\X$ be given. Define a new functor $S':\C\to\X$ as follows.
For an object $C$ of \C let $e_C$ be the composite of all the restriction
idempotents $\rst{Sf}$ where $f:C\to D$ is an arrow in \C with domain
$C$. Let $i_C:S'C\to SC$ and $r_C:SC\to S'C$ be the splitting of $e_C$.
Given an arrow $f:C\to D$, we have $e_D.Sf.e_C=Sf.e_C$, and so $Sf$
restricts to an arrow $S'f:S'C\to S'D$ satisfying $i_D.S'f=Sf.i_C$; and
this defines a functor $S':\C\to\X$ with a natural transformation
$i:S'\to S$. Since $S'$ lands in \TotX, we may form its limit
$p_C:L\to S'C$ in \TotX, and now $i_C p_C:L\to SC$ give a cone over
$S$ with total components; we shall show that it is a restriction limit
cone.

Let $q_C:M\to SC$ be the components of a lax cone over $S$, and write
$d:M\to M$ for the composite of the restriction idempotents $\rst{q_C}$. 
We must show that there is a unique arrow $f:M\to L$ satisfying
$i_C.p_C.f=q_C.d$.
Let $i:M'\to M$ and $r:M\to M'$ be a splitting of $d$. Each composite
$q_C i$ is total, and for an arrow $f:C\to D$ in \C we have
$q_D.i=Sf.q_C.\rst{q_D}i=Sf.q_C.i$; thus the $q_C i$ form the components
of a cone. For each $f:C\to D$ we have
$\rst{Sf}.q_C.i=q_C.i.\rst{Sf.q_C.i}=q_C.i.\rst{q_D.i}=q_C.i$,
and so $i_C.r_C.q_C.i=q_C.i$; but this means that $r_C.q_C.i$ is total,
and forms a cone over $S'$. Thus by the universal property of the
limit $p_C:L\to S'C$, there is a unique total map $g:M'\to L$ satisfying
$p_C.g=r_C.q_C.i$. Now
$i_C.p_C.gr=i_C.r_C.q_C.ir=q_C.ir=q_C.d$, so that $gr$ shows the
existence of an $f$. 

As for the uniqueness, let $f$ is any map
satisfying $i_C.p_C.f=q_C.d$; then $\rst{f}=\rst{i_C.p_C.f}=\rst{q_C.d}=
\rst{\rst{q_C}.d}=\rst{d}=d$. Now $fi$ satisfies
$i_C.p_C.fi=q_C.di=q_C.i$, and $\rst{fi}=\rst{\rst{f}i}=\rst{di}=\rst{i}=1$;
thus by the universal property of the limit $L$ in \TotX, we have 
$fi=g$, and now $f=f\rst{f}=fd=fir=gr$.
\endproof

Finally we observe that under a further assumption, restriction products
and splitting of restriction idempotents suffice to obtain all (finite)
limits in the restriction category.

Let \X be a cartesian restriction category. We say that an object $X$
is {\em separable} if the diagonal $\Delta:X\to X\t X$ is a restriction
monic; that is, if there is a map $r:X\t X\to X$ with $r\Delta=1$ and
$\Delta r=\rst{r}$. (Recall that such an $r$ is unique if it exists, and
is called the restriction retraction of $\Delta.$)

\begin{proposition}
If \X is a split cartesian restriction category in which every object is
separable then $\Total(\X)$ has all finite limits.
\end{proposition}

\proof
We already know that $\Total(\X)$ has finite products; it remains to
show that it has equalizers. Suppose then that $f,g:X\to Y$ are given
in $\Total(\X)$, let $h:X\to Y\t Y$ be the induced map, and $r:Y\t Y\to Y$
the restriction retraction of $\Delta:Y\to Y\t Y$. Now consider the
restriction idempotent $\rst{rh}$, and let $i:E\to X$ and $s:X\to E$
be its splitting. We shall show that $i$ is the desired equalizer. First
of all $i$ has a retraction, so is a monomorphism, and so in turn is total.
We must show that $fi=gi$, and that if $j$ is any total map with $fj=gj$
then $j$ factorizes through $i$.

Write $p,q:Y\t Y\to Y$ for the projections. Observe first that
$$p\rst{r}=p\Delta r=r=q\Delta r=q\rst{r}$$
and now
$$fi=fisi=f\rst{rh}i=ph\rst{rh}i=p\rst{r}hi=q\rst{r}hi=qh\rst{rh}i=gisi=gi.$$
On the other hand, if $j$ is total and $fj=gj$, then $hj=\Delta k$ for a 
(unique total) map $k$, and so
$$isj=\rst{rh}j=j\rst{rhj}=j\rst{r\Delta k}=j\rst{k}=j$$
and $sj$ gives the required factorization of $j$ through $i$.
\endproof


\section{Counital copy categories}
\label{sect:copy}


We have seen that there are many different ways of describing
the structure which we call a restriction category with restriction
products, but we shall actually add one more way to this list:
the {\em counital copy categories} which we introduce below.
(The slightly weaker structure of {\em copy category} will not
be considered in this paper.)


\subsection{Restriction products revisited}


The starting point is that if \X is a restriction category with
restriction products, then as a category, \X has a symmetric
monoidal structure, with tensor product given by restriction product.
The associativity isomorphism is the $\alpha$ appearing in the
definition of p-category, while the symmetry is the $\tau$. The 
unit is the restriction terminal object, and the unit constraint
$X\t T\cong X$ is the projection. In light of the coherence results
for monoidal categories \cite{CWM}, we shall allow ourselves to
omit explicit mention of the associativity isomorphisms, and write
as if the tensor product were strictly associative.

As observed by Carboni \cite{carboni}, for each object $X$ the diagonal map 
$\Delta:X\to X\t X$ is coassociative and cocommutative, and has a counit
given by $t_X:X\to T$. Thus every object has a canonical cocommutative
comonoid structure in the symmetric monoidal category. Furthermore,
since the $\Delta$ are natural, every morphism $f:X\to Y$ is a morphism
of cosemigroups, although it may not preserve the counit. It will preserve
the counit if it is a total map; conversely, if $f$ preserves the counit,
that is, if $t_Yf=t_X$, then $\rst{f}=\rst{t_Yf}=\rst{t_X}=1$, and so $f$
is total. Thus the total maps are precisely the counit-preserving ones.

There are two further further conditions which necessarily hold in 
a restriction category with restriction products: the diagonal 
$\Delta:T\to T\t T$ must be inverse to the unit isomorphism
$r=l:T\t T\to T$ of the monoidal structure, and the composite
$$\xymatrix{
X\t Y \ar[r]^-{\Delta\t\Delta} & X\t X\t Y\t Y \ar[r]^{1\t\tau\t1} &
X\t Y\t X\t Y}$$
must be $\Delta:X\t Y\to X\t Y\t X\t Y $.
Together these conditions say that $\Delta:X\to X\t X$ is not just a
natural transformation, but a {\em monoidal natural transformation}; it
can also be viewed as an instance of the ``middle four interchange'' law
for bicategories. To see that these conditions must hold in a restriction 
category with restriction products, it suffices to observe that they are 
all equations in \TotX, where the tensor product is a genuine product, and 
that such equations always hold in a symmetric monoidal category for which 
the tensor product is the categorical (cartesian) product. We call 
a symmetric monoidal category equipped with maps $\Delta:X\to X\t X$
which are monoidally natural, coassociative, cocommutative, and have
counits, a {\em counital copy category}. (In \cite{restiv} we shall have 
cause to look at a slightly weaker structure, called a {\em copy category}, 
in which the assumption that the cosemigroups have counits is dropped.)
It will turn out that a symmetric monoidal category can have at most one 
counital copy structure, and has such a structure if and only if it arises 
from a restriction category with restriction products.

If \V is an arbitrary symmetric monoidal category, let $\Copy(\V)$ be
the category whose objects are the cocommutative comonoids in \V, and
whose morphisms are the homomorphisms of cosemigroups. Then $\Copy(\V)$
has a canonical symmetric monoidal structure: the tensor product of
cocommutative comonoids $(C,\delta:C\to C\ot C,\epsilon:C\to I)$ and 
$(D,\delta:D\to D\ot D,\epsilon:D\to I)$ has underlying \V-object $C\ot D$,
with comultiplication and counit given by
$$\xymatrix{
C\ot D \ar[r]^-{\delta\ot\delta} & C\ot C\ot D\ot D \ar[r]^{1\ot\tau\ot1} &
C\ot D\ot C\ot D}$$
$$\xymatrix{
C\ot D \ar[r]^{\epsilon\ot\epsilon} & I\ot I \ar[r]^{r} & I. }$$
But now the map $\delta:C\to C\ot C$ in \V is a map 
$\Delta:(C,\delta,\epsilon)\to(C,\delta,\epsilon)\ot(C,\delta,\epsilon)$
in $\Copy(\V)$ which is coassociative, cocommutative, and counital by
definition of the objects of $\Copy(\V)$, natural by definition of morphisms 
in $\Copy(\V)$, and monoidally natural by definition of the monoidal
structure on $\Copy(\V)$. Thus $\Copy(\V)$ is a counital copy category.
On the other hand, there is an evident forgetful functor 
$U:\Copy(\V)\to\V$ which strictly preserves the symmetric monoidal structure,
and if \V\ {\em is} a counital copy category, then this $U$ is clearly an 
equivalence of categories. This proves:

\begin{proposition}
The counital copy categories are precisely the symmetric monoidal categories
of the form $\Copy(\V)$ for some symmetric monoidal \V.
\end{proposition}

We conclude:

\begin{theorem}
The following structures on a category \X are equivalent:
\begin{enumerate}[(i)]
\item restriction category with restriction products;
\item p-category with a one-element object;
\item partial cartesian category in the sense of Curien and Obtulowicz;
\item counital copy category;
\item symmetric monoidal structure with $U:\Copy(\X)\to\X$ an equivalence;
\item symmetric monoidal structure for which there exists some equivalence
$\X\simeq\Copy(\V)$.
\end{enumerate}
All the structure is determined by either the restriction category
structure or the symmetric monoidal structure.
\end{theorem}


\subsection{Classified restriction categories and equational lifting
categories}


In this brief section we revisit the analysis in \cite{restii}
of classified restriction categories, in particular its connection
with the {\em equational lifting monads} of \cite{bfs}.

Let \C be a symmetric monoidal category, with tensor product $\ot$,
unit $I$, and symmetry $\tau$. The associativity and unit isomorphisms
will be suppressed where possible. A {\em symmetric monoidal monad} \cite{kock}
on \C is a monad $T=(T,\eta,\mu)$ equipped with a natural transformation 
$\phi_{A,B}:TA\ot TB\to T(A\ot B)$ satisfying the equations 
$$\xymatrix{
TA\ot TB \ar[r]^{\phi_{A,B}} \ar[d]_{\tau} & T(A\ot B) \ar[d]^{\tau} &
A\ot B \ar[r]^{\eta_A\ot\eta_B} \ar[dr]_{\eta_{A\ot B}} & 
TA\ot TB \ar[d]^{\phi_{A,B}} \\
TB\ot TA \ar[r]_{\phi_{B,A}} & T(B\ot A) && T(A\ot B) }$$
$$\xymatrix{T^2A\ot T^2B \ar[r]^{\phi_{TA,TB}} \ar[d]_{\mu_A\ot\mu_B} &
T(TA\ot TB) \ar[r]^{T\phi_{A,B}} & T^2(A\ot B) \ar[d]^{\mu_{A\ot B}} \\
TA\ot TB \ar[rr]_{\phi_{A,B}} && T(A\ot B). }$$

In fact the structure on $T$ of symmetric monoidal monad can be
given either by $\phi$, or by a natural family of maps 
$\psi_{A,B}:A\t TB\to T(A\t B)$ satisfying equations given in \cite{kock}.
One obtains $\psi_{A,B}$ from
$\phi_{A,B}$ by composing with $\eta_A\t1_{TB}$, and one obtains $\phi_{A,B}$
from $\psi_{A,B}$ as the composite
$$\xymatrix{
TA\t TB \ar[r]^{\psi_{TA,B}} & T(TA\t B) \ar[r]^{T\tau} & 
T(B\t TA) \ar[r]^{T\psi_{B,A}} & T^2(B\t A) \ar[r]^{T^2\tau} &
T^2(A\t B) \ar[r]^{\mu_{A\t B}} & T(A\t B).}$$ 
When the maps $\psi$ are used rather than $\phi$, one sometimes speaks
of a {\em commutative strong monad} rather than a {\em symmetric monoidal
monad}.

A symmetric monoidal monad $T$ on \C induces a symmetric monoidal structure
on the Kleisli category $\C_T$. If we regard the objects
of $\C_T$ as being the objects of \C, and arrows in $\C_T$ from $A$ to $B$
as being arrows in \C from $A$ to $TB$, then the product of objects
$A$ and $A'$ is $A\ot A'$, while the product of arrows $f:A\to TB$ and
$f':A'\to TB'$ is the composite
$$\xymatrix{
A\ot A' \ar[r]^{f\ot f'} & TB\ot TB' \ar[r]^{\phi_{B,B'}} & T(B\ot B').}$$
The left adjoint $I:\C\to\C_T$ strictly preserves the symmetric monoidal
structure.

If \C is not just symmetric monoidal, but a counital copy category, then
applying $I:\C\to\C_T$ to the cocommutative comonoid structures on objects 
of \C, one obtains a canonical cocommutative structure on each object of
$\C_T$. If the resulting copy maps $A\to A\ot A$ in $\C_T$ are natural,
then they will certainly be monoidally natural, so that $\C_T$ will be
a counital copy category. As for the naturality, this amounts to commutativity
of the exterior of
$$\xymatrix{
A \ar[d]_{\Delta} \ar[rr]^f && TB \ar[dl]_{\Delta} \ar[d]^{T\Delta} \\
A\ot A \ar[r]_{f\ot f} & TB\ot TB \ar[r]_{\phi_{B,B}} & T(B\ot B)}$$
for every $f:A\to TB$ in \C. Now the quadrilateral commutes by naturality 
of the copy maps in \C, so the exterior will commute provided that the 
triangular region does so. We therefore define a symmetric monoidal
monad $T$ on a counital copy category \C to be a {\em copy monad} if
$\phi_{B,B}\Delta=T\Delta$ for all objects $B$.

\begin{proposition}
If $T$ is a copy monad on a counital copy category \C, then $\C_T$
is a counital copy category.
\end{proposition}

\begin{example}\label{ex:prod}
If \D is a distributive category, then the monad $+1$ on \D is
a symmetric monoidal monad, via the maps
$$\xymatrix{(A+1)\t(B+1) \ar[r]^-{\delta^{-1}} & 
A\t B+A+B+1 \ar[r]^-{A\t B+!} & A\t B+1.}$$
The fact that $+1$ is a copy monad amounts to commutativity of the exterior
of
$$\xymatrix{
(A+1)\t(A+1) \ar[r]^{\delta^{-1}} & A\t A+A+A+1 \ar[dr]^{A\t A+!} \\
A+1 \ar[u]^{\Delta} \ar[r]_{\Delta+1} & A\t A+1 \ar[u]^{i_{14}} \ar[r]_1 &
A\t A+1.}$$
It follows that \Dpar is a counital copy category, and so that
$\Total(K_r(\Dpar))$ has finite products.
\end{example}

In \cite{bfs}, a symmetric monoidal monad on a category with finite
products \C is called an {\em equational lifting monad} if 
$\psi_{A,B}:A\t TB\to T(A\t B)$ satisfies
$$\xymatrix{TA \ar[r]^{\Delta} \ar[d]_{T\Delta} & TA\t TA\ar[d]^{\psi_{TA,A}}\\
T(A\t A) \ar[r]_{T(\eta_A\t 1)} & T(TA\t A).}$$
We observe that this implies commutativity of
$$\xymatrix{
TA \ar[r]^{T\Delta} \ar[d]_{\Delta} \ar@/^2pc/[rr]^{T\Delta} &
T(A\t A) \ar[r]^{T\tau} \ar[d]^{T(\eta_A\t1)} & 
T(A\t A) \ar[d]_{T(1\t\eta_A)} \ar[dr]_{T\eta_{A\t A}} \ar[drr]^1 \\
TA\t TA \ar[r]^{\psi_{TA,A}} \ar@/_2pc/[rrrr]_{\phi_{A,A}} &
T(TA\t A) \ar[r]^{T\tau} & T(A\t TA) \ar[r]_{T\psi_{A,A}} &
T^2(A\t A) \ar[r]^{\mu_{A\t A}} & T(A\t A) }$$
which is to say that $T$ is a copy monad. Thus every equational lifting
monad is a copy monad.

\begin{question}
Is there a copy monad on a category with finite products which is not an
equational lifting monad? 
\end{question}

In \cite{restii}, we defined the notion of a {\em classifying monad}
on a category \C, and gave various characterizations. To give a monad
$T$ the structure of a classifying monad is to give its Kleisli category
$\C_T$ the structure of a restriction category for which $F_T:\C\to\C_T$
takes its values among the total maps, and the components of the counit
$\epsilon_T:F_TU_T\to 1$ are restriction retractions. We now prove:

\begin{proposition}
An equational lifting monad is a classifying monad.
\end{proposition}

\proof
We have already seen that $\C_T$ is a counital copy category, and so in 
particular a restriction category. The restriction of $f:A\to TB$ is given by 
$$\xymatrix{
A \ar[r]^-{\<1,f\>} & A\t TB \ar[r]^{\eta_A\t1} & TA\t TB \ar[r]^{\phi_{A,B}} &
T(A\t B) \ar[r]^-{Tp} & TA. }$$
By \cite[Proposition~3.15]{restii}, $T$ will be a classifying monad if
and only if the restriction of $\eta_B f:A\to TB$ is $\eta_A$, for each
$f:A\to B$ in \C; and the restriction of $1:TA\to TA$ is $T\eta_A::TA\to T^2A$.
The restriction of $\eta_B f:A\to TB$ is given by
$$\xymatrix{
A \ar[r]^-{\<1,\eta_Bf\>} & A\t TB \ar[r]^{\eta_A\t 1_TB} & 
TA\t TB \ar[r]^{\phi_{A,B}} & T(A\t B) \ar[r]^-{Tp} & TA }$$
and
$$T(p)\phi_{A,B}(\eta_A\t1_{TB})\<1,\eta_Bf\>=
T(p)\phi_{A,B}(\eta_A\t\eta_B)\<1,f\>=T(p)\eta_{A\t B}\<1,f\>=
\eta_A p\<1,f\>=\eta_A$$
as required. For the latter, the restriction of $1:TA\to TA$ is
$$\xymatrix{TA \ar[r]^-{\Delta} & TA\t TA \ar[r]^{\eta_{TA}\t1} &
T^2A\t TA \ar[r]^{\phi_{TA,A}} & T(TA\t A) \ar[r]^-{Tp} & T^2A}$$
and 
$$T(p)\phi_{TA,A}(\eta_{TA}\t1)\Delta=T(p)\psi_{A,A}\Delta=
T(p)T(\eta_A\t1)T(\Delta)=T\eta_A$$
as required.
\endproof


\subsection{Distributive copy categories}


A counital copy category \X is a restriction category with restriction
products; if \X also has restriction coproducts and the canonical
maps $\delta:A\t B+A\t C\to A\t(B+C)$ are invertible for all objects $A$,
$B$, and $C$ then we call \X a {\em distributive copy category}.

\begin{proposition}
For a  counital copy category \X with restriction coproducts, the
following are equivalent:
\begin{enumerate}[{\em ($i$)}]
\item \X is a distributive copy category;
\item $\Total(\X)$ is a distributive category;
\item $K_r(\X)$ is a distributive copy category;
\item $\Total(K_r(\X))$ is a distributive category.
\end{enumerate}
\end{proposition}

\proof
The equivalence of ($i$) and ($ii$) is immediate from the fact 
that restriction products and restriction coproducts in \X are 
products and coproducts in $\Total(\X)$; the equivalence of ($iii$)
and ($iv$) is a special case of this. The fact that ($iii$) implies
($i$) is trivial; it remains only to show that if the canonical
map $A\t B+A\t C\to A\t(B+C)$ is invertible for every object in \X,
then it is so for every object in $K_r(\X)$. This follows easily
from the fact that the objects of $K_r(\X)$ are retracts of the
objects of \X.
\endproof

Since in a distributive category the unique map $0\to A\t 0$
is invertible for any object $A$, the proposition implies that
the same is true for a distributive copy category.

Our main result about distributive copy categories is:

\begin{theorem}
If \X is a counital copy category with restriction coproducts,
then \X is an extensive restriction category if and only if 
it is a distributive copy category and has a restriction zero.
\end{theorem}

\proof
If \X is an extensive restriction category with restriction products
then it has a restriction zero by definition of extensivity for 
restriction categories; and $\Total(K_r(\X))$ is extensive with
finite products, thus distributive, so that \X is a distributive
copy category by the proposition.

Suppose conversely that \X is a distributive copy category with a
restriction zero. We must show that every map $f:C\to 1+1$ has
a decision. Let $h$ be the composite
$$\xymatrix{C \ar[r]^-{\Delta} & C\t C \ar[r]^-{C\t f} & 
C\t(1+1) \ar[r]^-{\delta^{-1}} & C+C.}$$
Then 
$\nabla h=p(C\t f)\Delta=p(C\t\rst{f})\Delta=p\Delta\rst{(C\t\rst{f})\Delta}=
p\Delta\rst{f}=\rst{f}$, giving one condition for $h$ to be an
$f$-decision. The second follows from commutativity of:
$$\xymatrix{
C \ar[r]^{\Delta} \ar[d]_{\Delta} & C\t C \ar[r]^{C\t f} \ar[d]_{\Delta\t C} &
C\t(1+1) \ar[r]^{\delta^{-1}} \ar[d]_{\Delta\t(1+1)} & 
C+C \ar[d]^{\Delta+\Delta} \\
C\t C \ar[r]_{C\t\Delta} \ar[dr]_{C\t f} & C\t C\t C \ar[r]_-{C\t C\t f} &
C\t C\t(1+1) \ar[r]^{\delta^{-1}} \ar[d]_{C\t f\t(1+1)} &
C\t C+C\t C \ar[d]^{C\t f+C\t f} \\
& C\t(1+1) \ar[r]_-{C\t\Delta} \ar[d]_{\delta^{-1}} &
C\t(1+1)\t(1+1) \ar[r]_-{\delta^{-1}} & 
C\t(1+1)+C\t(1+1) \ar[d]^{\delta^{-1}+\delta^{-1}} \\
& C+C \ar[rr]_{i+j} && C+C+C+C.}$$
\endproof

Thus we have the following examples of distributive copy categories:

\begin{example}
\begin{enumerate}[(i)]
\item For a distributive category \D we saw in Example~\ref{ex:prod}
that \Dpar is a counital copy category and in Example~\ref{ex:Dex}
that it is an extensive restriction category. By the theorem, then,
\Dpar is a distributive copy category.
\item If \V is a symmetric monoidal category then $\Copy(\V)$ is
a counital copy category. If \V also has coproducts, and the tensor
product distributes over the coproducts, then $\Copy(\V)$ is a 
distributive copy category.
\item The category \CRng of commutative rings can of course be regarded
as the category of commutative monoids in the monoidal category \Ab of
abelian groups. We can therefore regard \CRng\op as the category of
cocommutative comonoids in the monoidal category \Ab\op. Now the tensor
product in \Ab\op distributes over coproducts, and so $\Copy(\Ab\op)$
is a distributive copy category. An object of $\Copy(\Ab\op)$ is a
cocommutative comonoid in \Ab\op; that is, a commutative ring. In fact
$\Copy(\Ab\op)$ is just $\CRngp\op$, where \CRngp is the category whose
objects are the commutative rings and whose morphisms are the functions 
preserving $+$, $\times$, and $0$, but not necessarily preserving $1$. 
Thus $\CRngp\op$ is a distributive copy category. It is not hard to see
that idempotents split in $\CRngp\op$, and that the category of total maps
is just $\CRng\op$, and so we recover the well-known fact that $\CRng\op$ is
extensive.
\end{enumerate}
\end{example}


\subsection{The extensive completion of a distributive category}
\label{sect:Dext}


In this section we apply the results obtained above to give
a description of the extensive completion of a distributive 
category. There is a 2-category \Dist of distributive categories,
functors preserving finite products and coproducts, and natural
transformations; and there is a full sub-2-category \Extpr of
\Dist consisting of those distributive categories which are 
also extensive. The inclusion has a left biadjoint, and the
value at a distributive category \D of this left biadjoint is
what we mean by the extensive completion of the distributive
category \D. An explicit construction of the extensive completion 
was given in \cite{extcomp}; here we shall give an alternative,
more conceptual, description.

Given a distributive category \D we have seen that there is a monad $+1$ on \D
whose Kleisli category \Dpar has a restriction structure. We may
now split the restriction idempotents in \Dpar, and then take the
total maps in this new restriction category, to give a category
$\Total(K_r(\Dpar))$. 
The image of the left adjoint $\D\to\Dpar$ lands in $\Total(\Dpar)$,
and if we compose the resulting functor $I:\D\to\Total(\Dpar)$ with the map 
$\Total(J):\Total(\Dpar)\to\Total(K_r(\Dpar))$
induced by the inclusion $J:\Dpar\to K_r(\Dpar)$, we obtain a functor
$N:\D\to\Total(K_r(\Dpar))$. It turns out that $N:\D\to\Total(K_r(\Dpar))$
exhibits $\Total(K_r(\Dpar))$ as the extensive completion of \D, as we
shall see below.

We saw in Example~\ref{ex:Dex} that $\Total(K_r(\Dpar))$ is
extensive, and we saw in Example~\ref{ex:prod} that it has
finite products, and so lies in \Extpr. The inclusion 
$J:\Dpar\to K_r(\Dpar)$ preserves restriction products and
restriction coproducts, and so the induced map 
$\Total(J):\Total(\Dpar)\to\Total(K_r(\Dpar))$ preserves products
and coproducts. The left adjoint $\D\to\Dpar$ preserves coproducts,
and the inclusion $\Total(\Dpar)\to\Dpar$ preserves  and reflects
them, so that $I:\D\to\Total(\Dpar)$ preserves coproducts. On the
other hand the left adjoint $\D\to\Dpar$ sends products to restriction
products, and so $I:\D\to\Total(\Dpar)$ also preserves products. 
Thus $N:\D\to\Total(K_r(\Dpar))$ preserves products and coproducts,
and so is a morphism in \Dist. 

It remains to check the universal property. To do this, we use the
theory of {\em effective completions of classifying monads} developed 
in \cite[Section~5]{restii}. Recall that a restriction category \X is
{\em classified} if the inclusion $\Total(\X)\to\X$ has a 
right adjoint $R$ for which the components $\epsilon_A:RA\to A$
are restriction retractions. The induced comonad on \X is called
the {\em classifying comonad}. A monad $T$ on a category \C 
was defined in \cite{restii} to be a {\em classifying monad} if it 
is equipped with the requisite structure to make the Kleisli category 
$\C_T$ into a classified restriction category whose classifying
comonad is the comonad induced by the Kleisli adjunction.
The classifying monad $T$ is said to be {\em effective} if 
the restriction category $\C_T$ is split and the left adjoint
$F_T:\C\to\C_T$ exhibits \C as the category of total maps in
$\C_T$. In other words, a classifying monad is effective if
it is the partial map classifier for a category of partial
maps. Various characterizations of effective classifying monads
were given in \cite[Theorem~5.8]{restii}. 

Given a classifying monad $T$ on a category \C, the restriction
category $\C_T$ is classified; the split restriction category
$K_r(\C_T)$ need not be classified in general, although it will
be if $T$ is an {\em interpreted classifying monad} in the sense
of \cite[Section~4]{restii}. The precise details of this definition
are unimportant in the present context, but it {\em is} important
to know that the monad $+1$ on a distributive category \D is 
an interpreted classifying monad, as observed in \cite[Example~4.15]{restii}.
For a general classifying monad $T$ there is nonetheless a universal
way to obtain a split classified restriction category from $\C_T$:
it is obtained by splitting more idempotents than just the restriction
ones, and is denoted by $\Kcr(\C_T)$; see \cite[Section~3.3]{restii}.

Since $\Kcr(\C_T)$ is a split classified restriction category, the
induced monad on $\Total(\Kcr(\C_T))$ is an effective classifying 
monad. It is in fact the universal way of associating an effective
classifying monad to the classifying monad $T$, in a sense made
precise in \cite[Section~5]{restii}, and so is called the {\em effective
completion} of the classifying monad $T$. In the case where the monad
$T$ is interpreted --- such as the monad $+1$ on a distributive category
--- then the effective completion may be described more simply
as  $\Total(K_r(\C_T))$. We shall use the universal
property of the effective completion to show that $\Total(K_r(\Dpar))$
is the extensive completion of the distributive category \D.

Consider distributive categories \D and \E, equipped with the 
corresponding interpreted classifying monads $+1$. A morphism
of classifying monads (in the sense of \cite{restii}) from 
$(\D,+1)$ to $(\E,+1)$ consists of a functor $H:\D\to\E$ equipped
with a family of maps $\phi:H(A+1)\to HA+1$ natural in $A$ and 
rendering commutative the following diagrams:
$$\xymatrix{
HA \ar[r]^-{Hi_A} \ar[dr]_{i_{HA}} & H(A+1) \ar[d]^{\phi_A} &
H(A+1+1) \ar[r]^{\phi_{A+1}} \ar[d]_{H(A+\nabla)} & H(A+1)+1 \ar[r]^{\phi_A+1}&
HA+1+1 \ar[d]^{HA+\nabla} \\
& HA+1 &
H(A+1) \ar[rr]_{\phi_A} && HA+1}$$
$$\xymatrix{
HA \ar[r]^-{H\<A,f\>} \ar[d]_{\<HA,Hf\>} &
H(A\t(B+1)) \ar[r]^{H\delta^{-1}} & H(A\t B+A) \ar[r]^-{H(\pi_1+!)} &
H(A+1) \ar[d]^{\phi_A} \\
HA\t H(B+1) \ar[r]_{HA\t\phi_B} & HA\t(HB+1) \ar[r]_{\delta^{-1}} &
HA\t HB+HA \ar[r]_-{\pi_1+!} & HA+1}$$
for all objects $A$ and all morphisms $f:A\to B+1$. A straightforward
argument shows that the condition involving a morphism $f:A\to B+1$ holds
for all such $f$ if and only if it holds for all $f$ with $B=1$; that is,
for all $a:A\to 1+1$. The resulting diagram is:
$$\xymatrix{
HA \ar[r]^-{H\<A,a\>} \ar[d]_{\<HA,Ha\>} &
H(A\t(1+1)) \ar[r]^{H\delta^{-1}} & H(A+A) \ar[r]^-{H(A+!)} &
H(A+1) \ar[d]^{\phi_A} \\
HA\t H(1+1) \ar[r]_{HA\t\phi_1} & HA\t(H1+1) \ar[r]_{\delta^{-1}} &
HA\t H1+HA \ar[r]_-{\pi_1+!} & HA+1}$$

Given morphisms $(H,\phi)$ and $(K,\psi)$ of classifying monads, a
transformation from $(H,\phi)$ to $(K,\psi)$ consists of a natural
transformation $\alpha:H\to K$ rendering commutative 
$$\xymatrix{
H(A+1) \ar[r]^{\phi_A} \ar[d]_{\alpha_{A+1}} & HA+1 \ar[d]^{\alpha_A+1} \\
K(A+1) \ar[r]_{\psi_A} & KA+1}$$
for all $A$. There is now a 2-category \Distcl consisting of the
distributive categories, the morphisms of classifying monads, and
the transformations of these. It is a full sub-2-category of the
2-categories {\sf icMnd} and {\sf cMnd} defined in \cite{restii}. The
biadjunctions constructed in \cite[Section~5.2]{restii} now show 
that $N:\D\to\Total(K_r(\Dpar))$ has a canonical structure 
of morphism of classifying monads, and exhibits $\Total(K_r(\Dpar))$
as the bireflection of \D into the full sub-2-category of \Distcl
consisting of the extensive categories with finite products.

We now turn to an analysis of the notion of morphism of classifying
monads. Suppose, as above,  that \D and \E are distributive.

\begin{lemma}\label{lemma:distcl}
If $H:\D\to\E$ preserves finite coproducts then there is a unique
$\phi$ making $H$ into a morphism of classifying monads, namely
$HA+!:HA+H1\to HA+1$.
\end{lemma}

\proof
If $\phi_A:HA+H1\to HA+1$ makes $H$ into such a morphism then it
must have the form $\<\alpha_A | \beta_A\>$ where $\alpha_A:HA\to HA+1$ 
and $\beta_A:H1\to HA+1$ are both natural in $A$. Compatibility of 
$\alpha_A$ with the first injection gives $\alpha_A=i_{HA}$, while 
naturality of $\beta_A$ gives commutativity of
$$\xymatrix{
H1 \ar[r]^-{\beta_0} \ar[dr]_{\beta_A} & H0+1 \ar[d]^{H!+1} \\ & HA+1.}$$
But $H0$ is initial, so $\beta_0:H1\to H0+1$ is the unique map, and
$\phi_A$ must be $HA+!$ as claimed.

Conversely, we must show that $\phi_A=HA+!$ does satisfy the various 
conditions. It is clearly natural and satisfies the compatibility 
conditions with $i_A:A\to A+1$ and $A+\nabla:A+1+1\to A+1$, so we
need only show compatibility with maps $a:A\to 1+1$. Commutativity
of 
$$\xymatrix{
HA\t H1+HA\t H1 \ar[r]^{\delta} \ar[d]_{\pi_1+\pi_1} & 
HA\t(H1+H1) \ar[r]^{HA\t\theta} & HA\t H(1+1) \\
HA+HA \ar[r]_{\theta} & H(A+A) \ar[r]_{H\delta} & 
H(A\t(1+1)), \ar[u]_{\<H\pi_1,H\pi_2\>} }$$
wherein $\theta$ denotes the canonical isomorphisms expressing the
fact that $H$ preserves coproducts, gives commutativity of 
$$\xymatrix{
&& HA\t(H1+1) \ar[r]^{\delta^{-1}} & HA\t H1+HA \ar[ddr]^{\pi_1+!} \\
& HA\t H(1+1) \ar[ur]^{HA\t\phi_1} \ar[r]^{HA\t\theta^{-1}} &
HA\t(H1+H1) \ar[u]_{HA\t(H1+!)} \ar[r]_{\delta^{-1}} &
HA\t H1+HA\t H1 \ar[d]_{\pi_1+\pi_1} \ar[u]^{1+\pi_1} \ar[dr]^{\pi_1+!} \\
HA \ar[ur]^{\<HA,Ha\>} \ar[r]_-{H\<A,a\>} & 
H(A\t(1+1)) \ar[r]_{H\delta^{-1}} \ar[u]_{\<H\pi_1,H\pi_2\>} & 
H(A+A) \ar[r]^{\theta^{-1}} \ar[dr]_{H(A+!)} &
HA+HA \ar[r]^{HA+!} \ar[dr]^{HA+H!} & HA+1 \\
&&&H(A+1) \ar[r]_{\theta^{-1}} & HA+H1 \ar[u]_{HA+!} }$$
\endproof

\begin{corollary}
If $H,K:\D\to\E$ preserve finite coproducts, then any 
natural transformation $\alpha:H\to K$ is a transformation
of classifying monads.
\end{corollary}

On the other hand, the coproduct-preserving functors are
not the only morphisms of classifying monads:

\begin{example} \label{ex:constant}
If $X$ is any object of \E, then the constant functor $\Delta X:\D\to\E$
at $X$ becomes a morphism of classifying monads if we define 
$\phi_A:H(A+1)\to HA+1$ to be the injection $X\to X+1$ for any $A$.
\end{example}

We may now prove

\begin{theorem}
The functor $N:\D\to\Total(K_r(\Dpar))$ exhibits 
$\Total(K_r(\Dpar))$ as the extensive completion of the
distributive category \D.
\end{theorem}

\proof
We know from \cite{restii} that composition with $N$ induces, for
any extensive category \E with products, an equivalence between the 
category of morphisms of classifying monads from \D to \E and
the category of morphisms of classifying monads from $\Total(K_r(\Dpar))$
to \E. It remains to prove that if $G:\D\to\E$ preserves finite products
and coproducts, and $(H,\phi):\Total(K_r(\Dpar))\to\E$ is the induced 
morphism of  classifying monads, then $H$ preserves finite products 
and coproducts. But since $H$ may be constructed as the
composite of $\Total(K_r(G_{+1})):\Total(K_r(\Dpar))\to\Total(K_r(\E_{+1}))$
and the canonical equivalence $\Total(K_r(\E_{+1}))\to\E$, it will
suffice to prove that $\Total(K_r(G_{+1}))$ preserves finite products
and coproducts. Since $G$ preserves coproducts, $G_{+1}$ preserves
restriction coproducts, and so $\Total(K_r(G_{+1}))$ preserves coproducts 
by Proposition~\ref{prop:cocart}. Similarly $G_{+1}$ preserves 
restriction products and so $\Total(K_r(G_{+1}))$ preserves products 
by Proposition~\ref{prop:cart}.
\endproof




\begin{thebibliography}{99}


\bibitem{bfs}
A.~Bucalo, C.~F\"uhrmann, and A.~Simpson
\newblock Equational lifting monads, 
\newblock {\em Theoretical Computer Science}, to appear.

\bibitem{carboni} 
A.~Carboni,
\newblock Bicategories of partial maps, 
\newblock {\em Cah. de Top. Geom. Diff.}, 28:111--126, 1987.

\bibitem{matrices}
A.~Carboni,  Matrices, relations, and group representations
{\em J. Algebra} 136:497--529, 1991.

\bibitem{CLW}
A.~Carboni, Stephen Lack, and R.F.C.~Walters,
\newblock Introduction to extensive and distributive categories,
\newblock {\em J. Pure Appl. Algebra} 84:145--158, 1993.

\bibitem{cockett}
J.R.B. Cockett,
\newblock Introduction to distributive categories,
\newblock {\em Math. Structures Comput. Sci.} 3:277--307, 1993.

\bibitem{resti}
J.R.B. Cockett and Stephen Lack,
\newblock Restriction categories I: Categories of partial maps, 
\newblock {\em Theoretical Computer Science}, 270:223--259, 2002.

\bibitem{restii}
J.R.B. Cockett and Stephen Lack,
\newblock Restriction categories II: Partial map classification,
\newblock {\em Theoretical Computer Science}, 294:61--102, 2003.

\bibitem{restiv}
J.R.B. Cockett and Stephen Lack,
\newblock Restriction categories IV: Enriched restriction categories, in
preparation.

\bibitem{extcomp}
J.R.B. Cockett and Stephen Lack,
\newblock The extensive completion of a distributive category,
\newblock {\em Theory Appl. Categ.} 8:541--554, 2001.


\bibitem{curien}
P.-L. Curien and A. Obtulowicz,
\newblock Partiality, Cartesian closedness, and toposes,
\newblock {\em Inform. and Comput.} 80:50--95, 1989.

\bibitem{dipaola-heller} 
R.A. Di Paola and A. Heller, 
\newblock Dominical categories: recursion theory without elements, 
\newblock {\em J.~Symbolic Logic}, 52:595--635, 1987.

\bibitem{cats-alligators}
Peter J. Freyd and Andr\'e Scedrov,
{\em Categories, allegories} North-Holland Mathematical Library, 39, 
North-Holland Publishing Co., Amsterdam, 1990.




\bibitem{hoehnke}
H.-J. Hoehnke,
\newblock On partial algebras,
\newblock {\em Colloq. Math. Soc. J\'anos Bolyai} 29:373--412, 1977.




\bibitem{kock} 
A. Kock, 
\newblock Strong functors and monoidal monads, 
\newblock {\em Arch. Math.}, 23:113--120, 1972.


\bibitem{CWM}
Saunders Mac~Lane,
\newblock {\em Categories for the Working Mathematician},
\newblock Springer-Verlag, New York-Heidelberg-Berlin, 1971.





\bibitem{robinson-rosolini} 
E.P. Robinson and G. Rosolini, 
\newblock Categories of partial maps,
\newblock {\em Information and computation\/} 79:94--130, 1988.



\end{thebibliography}
\end{document}